\documentclass[preprint,12pt]{elsarticle}

\usepackage{epstopdf}
\usepackage{amsmath}
\usepackage{amsthm}
\usepackage{bm}
\usepackage{svg}
\usepackage{CJK}
\usepackage{enumerate}
\usepackage{color} 
\usepackage{amssymb}
\usepackage{algorithm}
\usepackage{algorithmic}

\usepackage{geometry}
\usepackage{amsfonts}
\usepackage{indentfirst}
\usepackage{mathtools}
\usepackage{setspace}
\usepackage{subfigure}
\usepackage{multirow}

\newtheorem{theorem}{Theorem}[section]

\newtheorem{prop}[theorem]{Proposition}

\newtheorem{remark}[theorem]{Remark}
\newtheorem{assumption}[theorem]{Assumption}

\usepackage{lineno}

\journal{TBD}

\begin{document}

\begin{frontmatter}

\title{Crank-Nicolson-type iterative decoupled algorithms for Biot's consolidation model using total pressure}

\author[sustech]{Huipeng Gu}
\author[morgan]{Mingchao Cai}
\author[sustech]{Jingzhi Li}

\address[sustech]{Department of Mathematics, Southern University of Science and Technology, Shenzhen, Guangdong 518055, China}
            
\address[morgan]{Department of Mathematics, Morgan State University, Baltimore, MD 21251, USA}

\begin{abstract}
In this work, we develop Crank-Nicolson-type iterative decoupled algorithms for a three-field formulation of Biot's consolidation model using total pressure. We begin by constructing an equivalent fully implicit coupled algorithm using the standard Crank-Nicolson method for the three-field formulation of Biot's model.
Employing an iterative decoupled scheme to decompose the resulting coupled system, we derive two distinctive forms of Crank-Nicolson-type iterative decoupled algorithms based on the order of temporal computation and iteration: a time-stepping iterative decoupled algorithm and a global-in-time iterative decoupled algorithm. 
Notably, the proposed global-in-time algorithm supports a partially parallel-in-time feature.
Capitalizing on the convergence properties of the iterative decoupled scheme, both algorithms exhibit second-order time accuracy and unconditional stability. 
Through numerical experiments, we validate theoretical predictions and demonstrate the effectiveness and efficiency of these novel approaches. 
\end{abstract}


\end{frontmatter}


\section{Introduction}
\label{Intro}

Biot's consolidation model, formulated by Biot in the mid-20th century \cite{biot1941general, biot1955theory}, elucidates the interplay between fluid flow and mechanical deformation within porous media. This model finds extensive applications in various fields such as geomechanics, geotechnical engineering, petroleum engineering, and biomechanics \cite{ju2020parameter}.
 
In this work, we investigate the three-field formulation of Biot's model \cite{oyarzua2016locking, lee2017parameter} in the following manner.
\begin{align}
    - \mbox{div} ( 2\mu \varepsilon (\bm{u})  -  \xi \mathbb{I} ) & = \bm{f} \quad \quad  \mbox{in} \ \Omega \times (0,T], \label{threefield1} 
    \\
    \mbox{div} \bm{u} + \frac{1}{\lambda} \xi - \frac{\alpha}{\lambda} p & = 0 \quad \quad \mbox{in} \ \Omega \times (0,T], \label{threefield2} 
    \\
    \left( c_0 + \frac{\alpha^2}{\lambda} \right) \partial_t p - \frac{\alpha}{\lambda} \partial_t \xi - \mbox{div} ( k_p \nabla p ) & = g \quad \quad \mbox{in} \ \Omega \times (0,T],  \label{threefield3} 
    \\
    \bm{u} = \bm{0} \quad \ \ \mbox{on} \ \Gamma_{\bm{u}} \times (0,T], \ \ \quad
     ( 2\mu \varepsilon (\bm{u})  -  \xi \mathbb{I} ) \bm{n} & = \bm{f}_1 \quad \ \ \mbox{on} \ \Gamma_{\sigma} \times (0,T],   \label{bc1} 
    \\
    p=0 \quad \quad \mbox{on} \ \Gamma_p \times (0,T], \quad \quad \quad \quad
    ( k_p \nabla p ) \cdot \bm{n} & = g_1 \quad \ \ \mbox{on} \ \Gamma_q \times (0,T], \label{bc4}
    \\
    \bm{u}(0) = \bm{u}^0, \quad p(0) = p^0, \quad \xi(0) = \alpha p^0 - \lambda & \mbox{div} \bm{u}^0  \quad \ \ \mbox{in} \ \Omega.  \label{ic}
\end{align}
Here, $\Omega$ is a bounded polygonal domain in $\mathbb{R}^d$ ($d= 2$ or $3$) with boundary $\partial\Omega = \Gamma_{\bm{u}} \cup \Gamma_{\sigma} = \Gamma_p \cup \Gamma_{q}$ with $\lvert \Gamma_{\bm{u}} \rvert>0, \lvert \Gamma_p \rvert >0$, $\Gamma_{\bm{u}} \cap \Gamma_{\sigma} = \Gamma_p \cap \Gamma_{q} = \emptyset$, and $T>0$ is the final time. Primary unknowns are the displacement vector of the solid $\bm{u}$, the fluid pressure $p$, and the total pressure $\xi = - \lambda \mbox{div} \bm{u} + \alpha p$. The term $\varepsilon(\bm{u}) = \frac{1}{2} [ \nabla \bm{u} + (\nabla \bm{u})^T ]$ denotes the strain tensor, $\mathbb{I}$ is the identity matrix, $\bm{f}$ is the body force, $g$ is a source or sink term, $c_0 \geq 0$ is the specific storage coefficient, $\alpha > 0$ is the Biot-Willis constant which is close to $1$, $k_p > 0$ represents the hydraulic conductivity, $\bm{n}$ is the unit outward normal to the boundary, and Lam{\'e} constants $\lambda$ and $\mu$ are computed from Young's modulus $E$ and the Poisson ratio $\nu$:
\begin{align*}
	\lambda = \frac{E\nu}{(1+\nu)(1-2\nu)},\ \ \ \mu = \frac{E}{2(1+\nu)}.
\end{align*}
In the above system, equations \eqref{threefield1} and \eqref{threefield2} describe the force equilibrium of the solid phase, and equation \eqref{threefield3} represents the mass conservation of the fluid phase.
The well-posedness analysis for the problem \eqref{threefield1}-\eqref{bc4} can be found in the works \cite{vzenivsek1984existence,showalter2000diffusion,oyarzua2016locking}.

In addressing the poroelasticity problem, three distinct numerical methods are employed in the literature: the fully coupled algorithm, the decoupled algorithm via time extrapolation, and the iterative decoupled algorithm. The fully coupled algorithm simultaneously solves the flow and mechanical equations as a unified system \cite{feng2018analysis, qi2021finite, burger2021virtual, gu2023priori}. 
Although it provides the advantages of unconditional stability and optimal convergence, it is notably hindered by its computational complexity. The decoupled algorithm tackles the flow and mechanical equations separately, in sequential time steps \cite{chaabane2018splitting, altmann2022decoupling, cai2023some}. This approach significantly enhances computational efficiency, although it does introduce a slight compromise in terms of stability. The iterative decoupled algorithm follows a different path by first decoupling the original problem and then iteratively addressing the subproblems \cite{kim2011stability, mikelic2013convergence, mikelic2014numerical,gu2023iterative, cai2023combination}. This strategy reduces computational demands while preserving optimality in both stability and accuracy. Among the various iterative decoupled methods, the fixed-stress splitting approach has garnered increasing attention due to its unconditionally stable and globally convergent properties \cite{mikelic2013convergence}. The iterative approach adopted in this study follows the fixed-stress splitting method. Notably, our algorithms are grounded in the three-field formulation (\ref{threefield1})-(\ref{threefield3}), departing from the conventional two-field formulation \cite{mikelic2013convergence, cai2023combination}. Utilizing the three-field formulation offers distinct advantages; specifically, classical Stokes inf-sup stable elements can be applied seamlessly, eliminating the need for additional stabilization \cite{lee2017parameter, oyarzua2016locking, gu2023iterative, cai2023some}.

As the poroelasticity problem inherently involves time-dependent aspects, achieving accuracy in the temporal domain is increasingly emphasized. Some higher-order methods in the coupled algorithm are well-established and documented in the literature \cite{wanner1996solving, kim2011stability, fu2019high, ge2019stabilized1}. However, employing higher-order methods in the decoupled algorithm may introduce stability concerns, as discussed in \cite{altmann2022semi}. In the context of the iterative decoupled method, the exploration of higher-order methods is discussed in \cite{bause2016iterative, bause2017space, kocher2018mixed}, which involves the introduction of continuous and discontinuous Galerkin space-time finite elements. 
Two prevalent fully discrete approaches, distinguished by the order of temporal computation and iteration, are the time-stepping type \cite{both2017robust, storvik2019optimization} and the global-in-time type \cite{borregales2019partially, ahmed2020adaptive}. The latter, which prioritizes subproblem iteration before handling time discretization across the entire domain, is particularly attractive for its compatibility with parallel computing. However, existing studies using the backward Euler method achieve only first-order convergence in time, posing challenges for long-time problems. To overcome this, the multirate fixed-stress method \cite{almani2016convergence} has been employed to reduce costs. However, recent analyses \cite{ge2023multirate} reveal it achieves only first-order convergence. Our contribution lies in developing both time-stepping and global-in-time iterative decoupled Crank-Nicolson type algorithms, leveraging the iterative strategy and decoupling techniques for second-order, optimally convergent solutions.

The paper unfolds as follows: Section \ref{Sec2} delves into the variational formulation, finite element discretization, and a synopsis of the reformulated fully coupled algorithm employing the Crank-Nicolson method. In Section \ref{Sec3}, leveraging the reformulation, we introduce several iterative decoupled algorithms along with their corresponding convergence analyses. The numerical experiments conducted to validate the theoretical results are elucidated in Section \ref{Sec5}. Lastly, Section \ref{conclusion} encapsulates the drawn conclusions from this study and delineates potential avenues for future research.

\section{The weak form and coupled numerical schemes}
\label{Sec2}

\subsection{Preliminaries}
We denote by $H^r(\Omega)$ the standard Sobolev space of functions whose distributional derivatives of order up to $r$ belong to $L^2(\Omega)$, equipped with the norm $\| \cdot \|_{H^r(\Omega)}$. For any $\Gamma \subset \partial \Omega$, we also define $H^r_{0,\Gamma} (\Omega) = \{ v \in H^r(\Omega), \ v=0 \ \text{on} \ \Gamma \}$ as the subspace of $H^r(\Omega)$. The inner products in $L^2(\Omega)$ and $L^2(\partial \Omega)$ are denoted by $( \cdot , \cdot )$ and $\langle \cdot , \cdot \rangle$, respectively. In this paper, we use $C$ to denote a generic positive constant independent of mesh sizes and use $x \lesssim y$ to denote $x \leq Cy$.

We first introduce the function spaces: 
\begin{align*}
    \bm{V} = \bm{H}^1_{0,\Gamma_{\bm{u}}} (\Omega), \quad W = L^2(\Omega), \quad M = H^1_{0,\Gamma_p} (\Omega).
\end{align*}
The weak formulation of equations \eqref{threefield1}-\eqref{bc4} reads as follows: find $(\bm{u}, \xi, p) \in \bm{V} \times W \times M$ satisfies the initial conditions \eqref{ic} such that for a.e. $t \in [0,T]$
\begin{align}
    a_1(\bm{u},\bm{v})-b(\bm{v},\xi) & = ( \bm{f},\bm{v} ) + \langle \bm{f}_1, \bm{v} \rangle_{\Gamma_{\sigma}}, \quad   \forall \bm{v} \in \bm{V}, \label{vf1} \\
    b(\bm{u},\phi) + a_2(\xi,\phi) - c(p,\phi) & =  0, \quad \quad  \quad \quad \quad \quad \quad \quad \forall \phi \in W, \label{vf2} \\
    a_3(\partial_t p,\psi)-c(\psi,\partial_t \xi)+d(p,\psi) & = ( g, \psi ) + \langle g_1 , \psi \rangle_{\Gamma_q}, \ \quad \forall \psi \in M. \label{vf3}
\end{align} 
Here, the bilinear forms are given by
\begin{align*}
    & a_1(\bm{u},\bm{v}) = 2\mu \int_\Omega  \varepsilon (\bm{u}) :  \varepsilon (\bm{v}),  \ \ \ b(\bm{v},\phi) =  \int_\Omega \phi \  \mbox{div}\bm{v}, \\
    & a_2(\xi,\phi) = \frac{1}{\lambda} \int_\Omega \xi \phi,  \quad \quad \quad \quad \quad \ c(p , \phi) =  \frac{\alpha}{\lambda} \int_\Omega p \phi, \\
    & a_3(p,\psi) = \left( c_0+\frac{\alpha^2}{\lambda} \right) \int_\Omega p \psi, \ \ \ d(p,\psi) =  k_p \int_\Omega \nabla p \cdot \nabla \psi.
\end{align*}
We note that the Korn's inequality \cite{korn1991nitsche} holds on $\bm{V}$, i.e., there exists a constant $C_K = C_K(\Omega,\Gamma_{\bm{u}}) > 0$ such that
\begin{align}
	\| \bm{v} \|_{H^1(\Omega)} \leq C_K \| \varepsilon(\bm{v}) \|_{L^2(\Omega)}, \quad \forall \bm{v} \in \bm{V}. \label{korncon}
\end{align}
Furthermore, the following inf-sup condition \cite{brenner1993nonconforming} holds: there exists a constant $\beta > 0$ depending only on $\Omega$ and $\Gamma_{\bm{u}}$ such that
\begin{align}
	\sup_{\bm{v}\in \bm{V}} \frac{  b(\bm{v},\phi) }{\| \bm{v} \|_{H^1(\Omega)}} \geq \beta \| \phi \|_{L^2(\Omega)}, \ \ \ \forall \phi \in W. \label{infsupcon}
\end{align}

In the subsequent sections, we make the following assumptions regarding the regularity of the solutions and right-hand terms:
\begin{assumption}
\label{assump} 
We assume that the exact solutions $(\bm{u},\xi,p)$ of \eqref{threefield1}-\eqref{bc4} satisfy 
\begin{align*}
    & \bm{u} \in C ([0,T]; \bm{H}^1_{0,\Gamma_{\bm{u}}} (\Omega) \cap \bm{W}^{1,\infty}(\Omega) ) \cap H^1([0,T];\bm{H}^{k+1}(\Omega)) \cap H^3([0,T];\bm{H}^1(\Omega)), 
    \\
    & \xi \in C([0,T]; L^{2}(\Omega) \cap W^{1,\infty}(\Omega)) \cap H^1([0,T];H^k(\Omega)) \cap H^3([0,T];L^2(\Omega)), \\
    & p \in C([0,T]; H^1_{0,\Gamma_p} (\Omega) \cap W^{1,\infty}(\Omega)) \cap H^1([0,T];H^{l+1}(\Omega)) \cap H^3([0,T];L^2(\Omega)),
\end{align*}
and the right-hand terms $\bm{f}$, $\bm{f}_1$, $g$ and $g_1$ are sufficiently smooth. Here, $k \geq 2$ and $l \geq 1$ are two integers determined by the finite element spaces to be used.
\end{assumption}

\subsection{Semi-discrete numerical scheme} Let $\mathcal{T}_h$ be a partition of the domain $\Omega$ into triangles in $\mathbb{R}^2$ or tetrahedras in $\mathbb{R}^3$, and let $h$ be the maximum diameter over all elements in the mesh. For the pair $(\bm{u}, \xi)$, we employ the Taylor-Hood element pair $(\bm{V}_h, W_h)$. For the pressure $p$, we adopt the Lagrange finite element space $M_h$. The finite element spaces are defined as follows:
\begin{align*} 
	& \bm{V}_{h} \coloneqq \{ \bm{v}_h \in \bm{H}^1_{0,\Gamma_{\bm{u}}} (\Omega) \cap {\bm{C} }^0(\bar{\Omega}); \ \bm{v}_h |_E  \in {\bm P}_{k}(E), ~\forall E \in \mathcal{T}_h \},
	\\
	& W_h \coloneqq \{ \phi_h \in L^2(\Omega) \cap C^0(\bar{\Omega}); \ \phi_h |_E  \in  P_{k-1}(E), ~\forall E \in \mathcal{T}_h \},
	\\
	& M_h \coloneqq \{ \psi_h \in H^1_{0,\Gamma_p} (\Omega) \cap C^0(\bar{\Omega}); \ \psi_h |_E  \in P_{l}(E), ~\forall E \in \mathcal{T}_h  \}.
\end{align*}
Here, the Taylor-Hood element pair $(\bm{V}_{h}, W_h)$ satisfies the discrete inf-sup condition: there exists a positive constant $\Tilde{\beta}$ independent of $h$ such that
\begin{align}
    \sup_{\bm{v}_h \in \bm{V}_h} \frac{  b(\bm{v}_h,\phi_h) }{\| \bm{v}_h \|_{H^1(\Omega)}} \geq \Tilde{\beta} \| \phi_h \|_{L^2(\Omega)}, \ \ \ \forall \phi_h \in W_h. \label{dinfsupcon}
\end{align}
Given the initial conditions $(\bm{u}_h(0), \xi_h(0), p_h(0))$, the semi-discrete variational formulation can be written as follows: for a.e. $t \in [0, T]$, find $(\bm{u}_h(t), \xi_h(t), p_h(t)) \in \bm{V}_h \times W_h \times M_h$ with the given initial conditions such that
\begin{align}
    a_1(\bm{u}_h,\bm{v}_h)-b(\bm{v}_h,\xi_h) & = ( \bm{f},\bm{v}_h ) + \langle \bm{f}_1, \bm{v}_h \rangle_{\Gamma_{\sigma}},  \quad \ \forall \bm{v}_h \in \bm{V}_h, \label{c1} \\
    b(\bm{u}_h,\phi_h) + a_2(\xi_h,\phi_h) - c(p_h,\phi_h) & =  0, \quad  \quad \quad \quad \quad \quad \quad \quad \quad \forall \phi_h \in W_h, \label{c2} \\
    a_3(\partial_t p_h,\psi_h)-c(\psi_h,\partial_t \xi_h)+d(p_h,\psi_h) & = ( g,\psi_h ) + \langle g_1 , \psi_h \rangle_{\Gamma_q}, \ \quad  \forall \psi_h \in M_h. \label{c3}
\end{align} 

\subsection{Crank-Nicolson-type  fully coupled algorithms} 
For the temporal discretization, we adopt an equidistant partition of the interval $0 = t_0 < t_1 < \cdots < t_{N} = T$ with a uniform time step denoted as $\Delta t$. For simplicity, we introduce the notations $\bm{u}^n = \bm{u}(t_n)$, $\xi^n = \xi(t_n)$, $p^n = p(t_n)$ to denote the exact solution, and $\bm{u}_h^n = \bm{u}_h(t_n)$, $\xi_h^n = \xi_h(t_n)$, $p_h^n = p_h(t_n)$ to represent the discrete solution.

We initiate our approach by employing the standard Crank-Nicolson method to solve problem \eqref{c1}-\eqref{c3}, which can be expressed as follows: for $1 \leq n \leq N$, 
given $(\bm{u}_h^{n-1},\xi_h^{n-1},p_h^{n-1}) \in \bm{V}_h \times W_h \times M_h$, find $(\bm{u}_h^{n},\xi_h^{n},p_h^{n}) \in \bm{V}_h \times W_h \times M_h$ such that for a.e. $t \in [0,T]$
\begin{align}
    & a_1\left( \frac{ \bm{u}_h^{n} + \bm{u}_h^{n-1} }{2},\bm{v}_h \right)-b\left(\bm{v}_h, \frac{ \xi_h^{n} + \xi_h^{n-1} }{2} \right) 
    \nonumber \\
    & \quad \quad \quad
    = \frac{1}{2} \left( \bm{f}^{n} + \bm{f}^{n-1} ,\bm{v}_h \right) + \frac{1}{2} \left\langle  \bm{f}_1^{n} + \bm{f}_1^{n-1} , \bm{v}_h \right\rangle_{\Gamma_{\sigma}} , \quad \quad \quad \quad \quad \quad \ \forall \bm{v}_h \in \bm{V}_h, \label{ck1}
    \\
    & b\left(\frac{ \bm{u}_h^{n} + \bm{u}_h^{n-1} }{2},\phi_h \right)+a_2 \left(\frac{ \xi_h^{n} + \xi_h^{n-1} }{2},\phi_h \right) - c \left(\frac{ p_h^{n} + p_h^{n-1} }{2},\phi_h \right) = 0 
    , \quad \ \forall \phi_h \in W_h, \label{ck2} \\
    & a_3 \left( \frac{p_h^{n} - p_h^{n-1}}{\Delta t}, \psi_h \right )-c\left( \psi_h , \frac{ \xi_h^{n}-\xi_h^{n-1} }{\Delta t} \right) + d\left( \frac{ p_h^{n} + p_h^{n-1} }{2},\psi_h \right)
    \nonumber \\
    & \quad \quad \quad  
    = \frac{1}{2} \left( g^{n} + g^{n-1} , \psi_h \right) +  \frac{1}{2} \left\langle g_1^{n}+g_1^{n-1} , \psi_h \right\rangle_{\Gamma_q}, \quad \quad \quad \quad \quad \quad \ \forall \psi_h \in M_h. \label{ck3}
\end{align}

Next, we reformulate the fully coupled system \eqref{ck1}-\eqref{ck3} to eliminate redundancy in the poroelasticity problem. Firstly, we note that the reformulation is predicated on a reasonable assumption, specifically, that the initial conditions adhere to the following expressions:
\begin{align}
    & a_1(\bm{u}_h^{0},\bm{v}_h)-b(\bm{v}_h,\xi_h^{0} ) = ( \bm{f}^{0},\bm{v}_h ) + \langle \bm{f}_1^{0} , \bm{v}_h \rangle_{\Gamma_{\sigma}} , \quad \ \forall \bm{v}_h \in \bm{V}_h, \label{ic_1} \\ 
    & b(\bm{u}_h^{0},\phi_h)+a_2(\xi_h^{0},\phi_h) - c(p_h^{0},\phi_h ) = 0 
    , \quad \quad \quad \quad \quad \forall \phi_h \in W_h. \label{ic_2}
\end{align}
Then, incremental progress can be achieved as follows:
For the time index $n=1$, by subtracting equations \eqref{ic_1} and \eqref{ic_2} from \eqref{ck1} and \eqref{ck2}, respectively, we can solve the alternative but equivalent system as follows:
\begin{align}
    & a_1(\bm{u}_h^{1},\bm{v}_h)-b(\bm{v}_h,\xi_h^{1} ) = ( \bm{f}^{1},\bm{v}_h ) + \langle \bm{f}_1^{1} , \bm{v}_h \rangle_{\Gamma_{\sigma}} , \quad \quad \quad \quad \quad  \forall \bm{v}_h \in \bm{V}_h,  \label{ic_3} \\ 
    & b(\bm{u}_h^{1},\phi_h)+a_2(\xi_h^{1},\phi_h) - c(p_h^{1},\phi_h ) = 0 
    , \quad \quad \quad \quad \quad \quad \quad \quad \ \forall \phi_h \in W_h, \label{ic_4} \\
    & a_3 \left( \frac{p_h^{1} - p_h^{0}}{\Delta t}, \psi_h \right )-c\left( \psi_h , \frac{ \xi_h^{1}-\xi_h^{0} }{\Delta t} \right) + d\left( \frac{ p_h^{1} + p_h^{0} }{2},\psi_h \right)
    \nonumber \\
    & \quad \quad \quad  
    = \frac{1}{2} \left( g^{1} + g^{0} , \psi_h \right) +  \frac{1}{2} \left\langle g_1^{1}+g_1^{0} , \psi_h \right\rangle_{\Gamma_q}, \quad \quad \quad \quad \quad \forall \psi_h \in M_h. \label{ic_5}
\end{align}
For the time index $n=2$, once again, we subtract equations \eqref{ic_3} and \eqref{ic_4} from \eqref{ck1} and \eqref{ck2}, respectively. Then, the alternative system reads as:
\begin{align}
    & a_1(\bm{u}_h^{2},\bm{v}_h)-b(\bm{v}_h,\xi_h^{2} ) = ( \bm{f}^{2},\bm{v}_h ) + \langle \bm{f}_1^{2} , \bm{v}_h \rangle_{\Gamma_{\sigma}}, \quad \quad \quad \quad \quad  \forall \bm{v}_h \in \bm{V}_h,  \label{ic_6} \\ 
    & b(\bm{u}_h^{2},\phi_h)+a_2(\xi_h^{2},\phi_h) - c(p_h^{2},\phi_h ) = 0 
    , \quad \quad \quad \quad \quad \quad \quad \quad \ \forall \phi_h \in W_h, \label{ic_7} \\
    & a_3 \left( \frac{p_h^{2} - p_h^{1}}{\Delta t}, \psi_h \right )-c\left( \psi_h , \frac{ \xi_h^{2}-\xi_h^{1} }{\Delta t} \right) + d\left( \frac{ p_h^{2} + p_h^{1} }{2},\psi_h \right)
    \nonumber \\
    & \quad \quad \quad  
    = \frac{1}{2} \left( g^{2} + g^{1} , \psi_h \right) +  \frac{1}{2} \left\langle g_1^{2} +g_1^{1} , \psi_h \right\rangle_{\Gamma_q}, \ \quad \quad \quad \quad \forall \psi_h \in M_h. \label{ic_8}
\end{align}
Iterating through the above-described process for all time steps, we obtain the reformulated Crank-Nicolson coupled algorithm, as summarized in \textbf{Algorithm \ref{algo:Coupled_CN}}.

\begin{algorithm}
\raggedright
  \caption{: A reformulated Crank-Nicolson coupled algorithm}
  \label{algo:Coupled_CN}
  \textbf{Input:} initial information $(\bm{u}_h^{0},\xi_h^{0},p_h^{0}) \in \bm{V}_h \times W_h \times M_h$.
  \\
  \textbf{Output:} solution at the final time $(\bm{u}_h^{N},\xi_h^{N},p_h^{N}) \in \bm{V}_h \times W_h \times M_h$.
  \\
  \textbf{for} $n$ from $1$ to $N$
  \\
  \quad \quad find $(\bm{u}_h^{n},\xi_h^{n},p_h^{n}) \in \bm{V}_h \times W_h \times M_h$ such that
  \begin{align}
    & a_1(\bm{u}_h^{n},\bm{v}_h)-b(\bm{v}_h,\xi_h^{n} ) = ( \bm{f}^{n},\bm{v}_h ) + \langle \bm{f}_1^{n} , \bm{v}_h \rangle_{\Gamma_{\sigma}} , \quad \quad \quad \quad  \quad \ \forall \bm{v}_h \in \bm{V}_h, \label{d1}
    \\
    & b(\bm{u}_h^{n},\phi_h)+a_2(\xi_h^{n},\phi_h) - c(p_h^{n},\phi_h ) = 0
    , \quad \quad \quad \quad \quad \quad \quad \quad \quad  \forall \phi_h \in W_h, \label{d2} 
    \\
    & a_3 \left( \frac{p_h^{n} - p_h^{n-1}}{\Delta t}, \psi_h \right )-c\left( \psi_h , \frac{ \xi_h^{n}-\xi_h^{n-1} }{\Delta t} \right) + d \left( \frac{p_h^{n}+p_h^{n-1}}{2},\psi_h \right)
    \nonumber \\
    & \quad \quad \quad  = \frac{1}{2} ( g^{n} + g^{n-1} , \psi_h ) + \frac{1}{2} \langle  g_1^{n} + g_1^{n-1}  , \psi_h \rangle_{\Gamma_q}, \quad \quad \quad \quad \forall \psi_h \in M_h. \label{d3}
  \end{align} \\
  \textbf{end for}
\end{algorithm}

\vspace{1cm}

In contrast to the coupled system given by equations \eqref{ck1}-\eqref{ck3}, the reformulated Crank-Nicolson coupled algorithm omits the reliance on information from the preceding time step as seen in equations \eqref{d1} and \eqref{d2}. This simplification proves advantageous for the implementation of the iterative decoupled approach, particularly in the context of the global-in-time algorithm. 
The following theorem establishes an energy estimate for the solution $(\bm{u}_h^{N},\xi_h^{N},p_h^{N})$ of \textbf{Algorithm \ref{algo:Coupled_CN}}. Interested readers can refer to our prior work for the corresponding proof (see Thm 4.6 in \cite{gu2023priori}).

\begin{theorem}
\label{CoThm}
Let $(\bm{u}^N,\xi^N,p^N)$ and $(\bm{u}^N_h,\xi^N_h,p^N_h)$ be solutions of problem \eqref{c1}-\eqref{c3} and problem \eqref{d1}-\eqref{d3} at the final time $t_N = T$, respectively. Under \textbf{Assumption \ref{assump}}, there holds
\begin{align}
    & \| \bm{u}_h^N - \bm{u}^N \|_{H^1(\Omega)} + \| \xi_h^N - \xi^N \|_{L^2(\Omega)} 
    \nonumber \\
    & \lesssim (\Delta t)^2 \int_{0}^{T} \left( \|\partial_{ttt}\bm{u}\|_{H^1(\Omega)} + \|\partial_{ttt}\xi \|_{L^2(\Omega)} + \|\partial_{ttt}p \|_{L^2(\Omega)} \right) ds
    \nonumber \\
    & \quad  + h^{k} \int_{0}^{T}\left( \|\partial_{t} \bm{u} \|_{H^{k+1}(\Omega)} + \|\partial_{t} \xi \|_{H^{k}(\Omega)}  \right)ds + h^{l+1} \int_{0}^{T}  \|\partial_{t} p \|_{H^{l+1}(\Omega)} ds, \label{coupled1}
    \\
    & \| p_h^N - p^N \|_{H^1(\Omega)} 
    \lesssim (\Delta t)^2 \int_{0}^{T} \left( \|\partial_{ttt}\bm{u}\|_{H^1(\Omega)} + \|\partial_{ttt}\xi \|_{L^2(\Omega)} + \|\partial_{ttt}p \|_{L^2(\Omega)} \right) ds
    \nonumber \\
    & \quad + h^k \int_{0}^{T}\left( \|\partial_{t} \bm{u} \|_{H^{k+1}(\Omega)} + \|\partial_{t} \xi \|_{H^{k}(\Omega)}  \right)ds + h^l \int_{0}^{T} \|\partial_{t} p \|_{H^{l+1}(\Omega)}ds . \label{coupled2}
\end{align}
\end{theorem}

The reformulated Crank-Nicolson coupled algorithm demonstrates unconditional stability and optimal convergence. However, the direct solution of a three-by-three system involves considerable computational expense. Consequently, we are motivated to investigate decoupling strategies. In the subsequent section, we present several decoupled schemes aimed at mitigating the computational costs associated with the fully coupled system.

\section{Crank-Nicolson-type iterative decoupled algorithms}
\label{Sec3}

This section aims to develop several iterative decoupled algorithms utilizing Crank-Nicolson-type time discretization. Our objective is to establish relationships between the approximations generated by these iterative decoupled algorithms and those obtained through the fully coupled algorithm. By leveraging Theorem \ref{CoThm}, we connect the approximations from the iterative decoupled methods to the exact solution. The core concept is illustrated 
\begin{figure}
\centering
\includegraphics[width=1.0\textwidth,trim=0 750 0 10,clip]{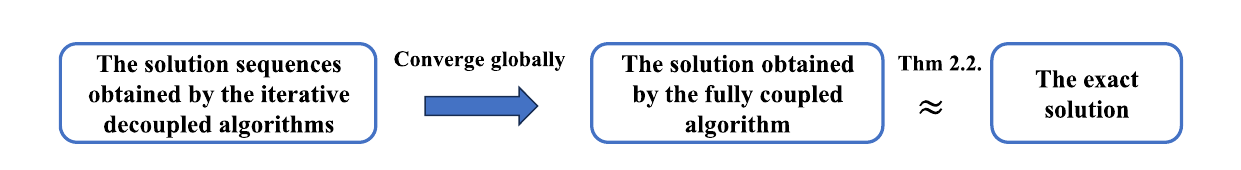}
\caption{
The connections among solution sequences produced by the iterative decoupled algorithms, the numerical solution derived from the reformulated Crank-Nicolson coupled algorithm, and the exact solution.}
\label{fig_ds}
\end{figure}
in Figure \ref{fig_ds}. Given the unconditional stability and second-order time accuracy of the reformulated Crank-Nicolson coupled algorithm, we infer that the Crank-Nicolson-type iterative decoupled algorithms also exhibit second-order time accuracy and stability. 

\subsection{A semi-discrete iterative decoupled algorithm}
To lay the foundation for our Crank-Nicolson-type iterative decoupled algorithms, we begin by introducing a semi-discrete iterative approach that decomposes the original problem into two subproblems: a reaction-diffusion problem and a generalized Stokes problem. The iterative process commences with an initial guess $\xi^0_h(t) \in W_h$. The resulting sequence $\{(\bm{u}_h^{i}(t),\xi_h^{i}(t),p_h^{i}(t))\}_{i \geq 0}$ is utilized to approximate the solution $(\bm{u}_h(t), \xi_h(t), p_h(t))$ governed by \eqref{c1}-\eqref{c3}. The $i$-th iteration unfolds as follows:
\\

\noindent \textbf{Step 1:} (The reaction-diffusion problem) Given $\xi^{i-1}_h(t) \in W_h$, for a.e. $t \in [0, T]$, find $p_h^i(t) \in M_h$ with $p_h^i(0) = p_h^0$ such that
\begin{align}
    a_3(\partial_t p_h^i,\psi_h) + d(p_h^{i},\psi_h) & = c(\psi_h,\partial_t \xi_h^{i-1}) + ( g,\psi_h ) + \langle g_1 , \psi_h \rangle_{\Gamma_q}, \ \ \forall \psi_h \in M_h. \label{fsw3}
\end{align}

\noindent \textbf{Step 2:} (The generalized Stokes problem) Given $p_h^i(t) \in M_h$, for a.e. $t \in [0,T]$, find $(\bm{u}_h^{i}(t),\xi_h^{i}(t)) \in \bm{V}_h \times W_h$ with $\bm{u}_h^i(0) = \bm{u}_h^0$, $\xi_h^i(0) = \xi_h^0$ such that
\begin{align}
    a_1(\bm{u}_h^{i},\bm{v}_h)-b(\bm{v}_h,\xi_h^{i}) & =  ( \bm{f},\bm{v}_h ) + \langle \bm{f}_1, \bm{v}_h \rangle_{\Gamma_{\sigma}}, \quad \quad \quad \quad \quad \forall \bm{v}_h \in \bm{V}_h, \label{fsw1} \\
    b(\bm{u}_h^{i},\phi_h) + a_2(\xi_h^{i},\phi_h) & = c(p_h^{i},\phi_h) , \ \quad \quad \quad \quad \quad \quad \quad \quad \quad \forall \phi_h \in W_h. \label{fsw2}
\end{align}

Next, we present the following theorem to demonstrate the convergence of the algorithm \eqref{fsw3}-\eqref{fsw2}. Importantly, based on the proposed three-field Biot'S model, we emphasize that the corresponding space-time convergence holds even when the storage coefficient $c_0$ is zero.

\begin{theorem}
\label{impthm}
The solution sequence $\{(\bm{u}_h^{i},\xi_h^{i},p_h^{i})\}_{i \geq 0}$ generated by equations \eqref{fsw3}-\eqref{fsw2} converges globally to $(\bm{u}_h,\xi_h,p_h)$, the solution of \eqref{c1}-\eqref{c3}.
We denote the iteration errors by $e_{\bm{u}}^i = \bm{u}_h^i - \bm{u}_h$, $e_{\xi}^i = \xi_h^i - \xi_h$, $e_{p}^i = p_h^i - p_h$. There holds
\begin{equation}
    \left[ \int_0^T \| \partial_{t} e_{\xi}^{i} \|_{L^2(\Omega)}^2 dt \right]^{\frac{1}{2}} \leq L \left[ \int_0^T \| \partial_{t} e_{\xi}^{i-1} \|_{L^2(\Omega)}^2  dt \right]^{\frac{1}{2}}, \label{thm1eq0}
\end{equation}
where $L$ is a positive constant strictly smaller than $1$ given as follows.
    \begin{align*}
        L = \frac{1}{ \left( \frac{c_0 \lambda}{\alpha^2} + 1  \right) \times \left( C(\beta,\mu) \lambda + 1 \right) }.
    \end{align*}
\end{theorem}

\begin{proof}
Subtracting the equation \eqref{c3} from the equation \eqref{fsw3} and setting the test function $\psi_h = \partial_t e_p^i$, we take the integral over $0$ up to $T$ and get
\begin{align}
    \int_0^T a_3(\partial_t e_p^i, \partial_t e_p^i) dt + \int_0^T d(e_p^{i},\partial_t e_p^i) dt = \int_0^T c(\partial_t e_p^i,\partial_t e_\xi^{i-1}) dt. \label{impthm:eq1}
\end{align}
Since $e_p^i(0) = 0$, we can use the Cauchy-Schwarz inequality to reformulate \eqref{impthm:eq1} as follows.
\begin{align*}
    & \left(c_0 + \frac{\alpha^2}{\lambda} \right) \int_0^T \| \partial_t e_p^i \|_{L^2(\Omega)}^2 dt + \frac{k_p}{2} \| \nabla e_{p}^{i} (T) \|_{L^2(\Omega)}^2 
    \nonumber \\
    & \leq \frac{1}{2} \left(c_0 + \frac{\alpha^2}{\lambda} \right) \int_0^T \| \partial_t e_p^i \|_{L^2(\Omega)}^2 dt + \frac{1}{2 \left(c_0 + \frac{\alpha^2}{\lambda} \right) } \frac{\alpha^2}{\lambda^2}  \int_0^T \| \partial_t e_\xi^{i-1} \|_{L^2(\Omega)}^2 dt,
\end{align*}
which leads directly to 
\begin{align}
    & \int_0^T \| \partial_t e_p^i \|_{L^2(\Omega)}^2 dt + \frac{k_p}{ \left(c_0 + \frac{\alpha^2}{\lambda} \right) } \| \nabla e_{p}^{i} (T) \|_{L^2(\Omega)}^2 \leq \frac{ \frac{\alpha^2}{\lambda^2} }{ \left(c_0 + \frac{\alpha^2}{\lambda} \right)^2 }  \int_0^T \| \partial_t e_\xi^{i-1} \|_{L^2(\Omega)}^2 dt. \label{impthm:eq2}
\end{align}

Subtracting the equations \eqref{c1} and \eqref{c2} from the equations \eqref{fsw1} and \eqref{fsw2} gives
\begin{align}
    a_1(\partial_t e_{\bm{u}}^{i},\bm{v}_h) - b(\bm{v}_h, \partial_t e_\xi^{i}) & =  0, \label{impthm:eq3} \\
    b(\partial_t e_{\bm{u}}^{i},\phi_h) + a_2(\partial_t e_\xi^{i} , \phi_h) & = c(\partial_t e_p^{i},\phi_h). \label{impthm:eq4}
\end{align}
Setting $\bm{v}_h = \partial_t e_{\bm{u}}^i$ in \eqref{impthm:eq3} and $\phi_h = \partial_t e_{\xi}^i$ in \eqref{impthm:eq4}, we take the integral over $0$ up to $T$ and then sum the resulting equations to get
\begin{align}
    \int_0^T a_1(\partial_t e_{\bm{u}}^{i},\partial_t e_{\bm{u}}^{i}) dt + \int_0^T a_2(\partial_t e_\xi^{i} , \partial_t e_\xi^{i} ) dt = \int_0^T c(\partial_t e_p^i,\partial_t e_\xi^{i}) dt.
\end{align}
We then apply the Cauchy-Schwarz inequality to derive that
\begin{align}
    & 2\mu \int_0^T \| \varepsilon(\partial_t e_{\bm{u}}^{i}) \|_{L^2(\Omega)}^2 dt + \frac{1}{\lambda} \int_0^T \| \partial_t e_{\xi}^{i} \|_{L^2(\Omega)}^2 dt
    \nonumber \\
    & \leq \frac{\alpha}{\lambda} \left( \int_0^T \| \partial_t e_p^i \|_{L^2(\Omega)}^2 dt \right)^{\frac{1}{2}} \left( \int_0^T \| \partial_t e_\xi^{i} \|_{L^2(\Omega)}^2 dt \right)^{\frac{1}{2}}. \label{impthm:ceq2}
\end{align}
One may apply the inf-sup condition \eqref{infsupcon} to derive the following inequality from \eqref{impthm:eq3}
\begin{align*}
    \beta \| \partial_t e_{\xi}^{i} \|_{L^2(\Omega)} \leq
    \sup_{\bm{v}\in \bm{V}} \frac{  b(\bm{v}, \partial_t e_{\xi}^{i} ) }{\| \bm{v} \|_{H^1(\Omega)}} = \sup_{\bm{v}\in \bm{V}} \frac{  a_1(\partial_t e_{\bm{u}}^{i},\bm{v}_h) }{\| \bm{v} \|_{H^1(\Omega)}} \leq 2 \mu C \| \varepsilon(\partial_t e_{\bm{u}}^{i}) \|_{L^2(\Omega)},
\end{align*}
which leads to the following inequality
\begin{align}
    C(\beta,\mu) \int_0^T \| \partial_t e_{\xi}^{i} \|_{L^2(\Omega)}^2 dt \leq 2\mu \int_0^T \| \varepsilon(\partial_t e_{\bm{u}}^{i}) \|_{L^2(\Omega)}^2 dt. \label{impthm:ceq3}
\end{align}

Finally, one can combine the equations \eqref{impthm:eq2}, \eqref{impthm:ceq2} and \eqref{impthm:ceq3} to derive
\begin{align*}
    & \left( C(\beta,\mu) \lambda + 1 \right)^2 \int_0^T \| \partial_t e_{\xi}^{i} \|_{L^2(\Omega)}^2 dt + \frac{k_p \alpha^2}{ \left(c_0 + \frac{\alpha^2}{\lambda} \right) } \| \nabla e_{p}^{i} (T) \|_{L^2(\Omega)}^2 
    \\ & 
    \leq \frac{ \frac{\alpha^4}{\lambda^2} }{ \left(c_0 + \frac{\alpha^2}{\lambda} \right)^2 }  \int_0^T \| \partial_t e_\xi^{i-1} \|_{L^2(\Omega)}^2 dt,
\end{align*}
leading to the following contraction
\begin{align}
    \left[ \int_0^T \| \partial_t e_{\xi}^{i} \|_{L^2(\Omega)}^2 dt \right]^{\frac{1}{2}} \leq \frac{1}{ \left( \frac{c_0 \lambda}{\alpha^2} + 1  \right) \times \left( C(\beta,\mu) \lambda + 1 \right) } \left[ \int_0^T \| \partial_t e_{\xi}^{i-1} \|_{L^2(\Omega)}^2 dt \right]^{\frac{1}{2}}.
\end{align}
Returning to \eqref{impthm:eq2}, we conclude that the error term $\| \nabla e_{p}^{i} (T) \|_{L^2(\Omega)}^2$ vanishes as $i \to \infty$. Given the uniqueness and existence of the solution, we conclude that the sequence $\{(\bm{u}_h^{i},\xi_h^{i},p_h^{i})\}_{i \geq 0}$ converges globally to $(\bm{u}_h, \xi_h, p_h)$.
\end{proof}

\subsection{A time-stepping iterative decoupled algorithm}

\begin{figure}[h]
\centering
\includegraphics[width=0.95\textwidth,trim=0 590 0 0,clip]{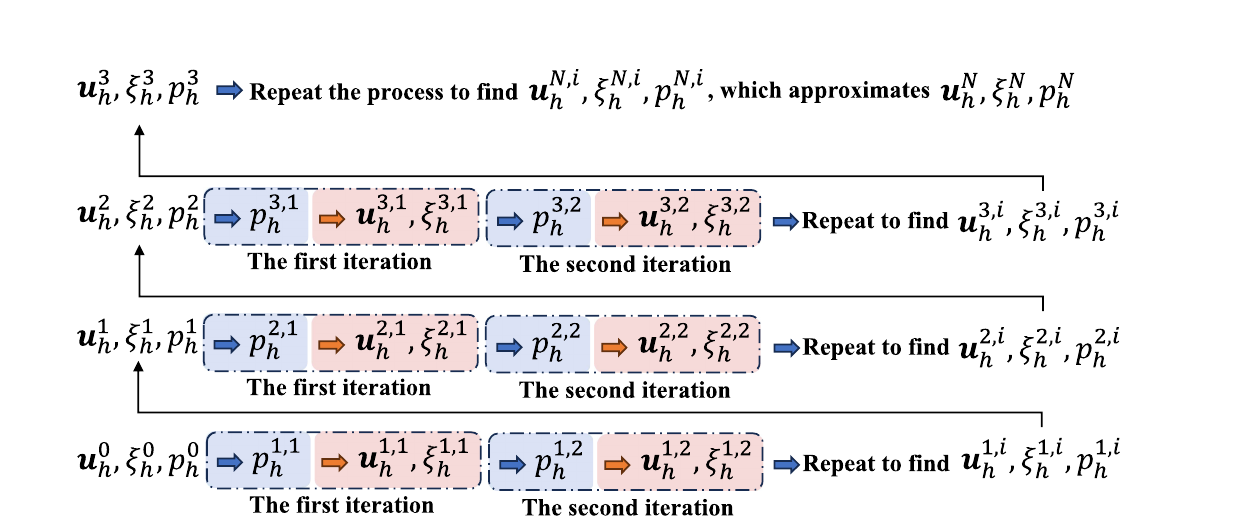}
\caption{
Progression in the time-stepping iterative decoupled algorithm.}
\label{fig_ts}
\end{figure}

The time-stepping method, also known as the time-marching method, stands as a widely adopted approach for addressing time-dependent differential equations. This method involves sequential problem-solving within distinct subintervals, such as $[t_0, t_1]$, $[t_1, t_2]$, $\cdots$, $[t_{N-1}, t_N]$. In \textbf{Algorithm \ref{algo:TIDA_CN}}, we deploy the Crank-Nicolson iterative decoupled algorithm employing the time-stepping method. The corresponding flowchart is depicted in Figure \ref{fig_ts}.
Following the idea in \cite{gu2023iterative}, we provide a convergence analysis for \textbf{Algorithm \ref{algo:TIDA_CN}} in Theorem \ref{thm2}.
Let $(\bm{u}_h^n, \xi_h^n, p_h^n)$ and $(\bm{u}_h^{n,i}, \xi_h^{n,i}, p_h^{n,i})$ represent the solutions of problem \eqref{d1}-\eqref{d3} and problem \eqref{ts3}-\eqref{ts2}, respectively. We also define the $i$-th iteration errors as $e_{\bm{u},s}^{n,i} = \bm{u}_h^{n,i} - \bm{u}_h^n$, $e_{\xi,s}^{n,i} = \xi_h^{n,i} - \xi_h^n$, $e_{p,s}^{n,i} = p_h^{n,i} - p_h^n$.

\begin{algorithm}[ht]
\raggedright
  \caption{: A time-stepping Crank-Nicolson iterative decoupled algorithm}
  \label{algo:TIDA_CN}
  \textbf{Input:} initial information $(\bm{u}_h^{0},\xi_h^{0},p_h^{0}) \in \bm{V}_h \times W_h \times M_h$.
  \\
  \textbf{Output:} solution at the final time $(\bm{u}_h^{N},\xi_h^{N},p_h^{N}) \in \bm{V}_h \times W_h \times M_h$.
  \\
  \textbf{for} $n$ from $1$ to $N$
  \\
  \quad \quad set $i=0$, $(\bm{u}_h^{n,0},\xi_h^{n,0},p_h^{n,0})=(\bm{u}_h^{n-1},\xi_h^{n-1},p_h^{n-1})$.
  \\
  \quad \quad \textbf{do} 
  \\
  \quad \quad \quad \quad set $i = i + 1$.
  \\
  \quad \quad \quad \quad \textbf{Step a:} find $p_h^{n,i}\in M_h$ such that
  \begin{align}
    & a_3 \left( \frac{p_h^{n,i} - p_h^{n-1}}{\Delta t}, \psi_h \right )-c\left( \psi_h , \frac{ \xi_h^{n,i-1}-\xi_h^{n-1} }{\Delta t} \right) + d \left( \frac{p_h^{n,i}+p_h^{n-1}}{2},\psi_h \right)
    \nonumber \\
    & \quad \quad = \frac{1}{2} ( g^{n} + g^{n-1} , \psi_h ) + \frac{1}{2} \langle  g_1^{n} + g_1^{n-1}  , \psi_h \rangle_{\Gamma_q}, \quad \quad \quad \ \forall \psi_h \in M_h. \label{ts3}
  \end{align}
  \\
  \quad \quad \quad \quad \textbf{Step b:} find $(\bm{u}_h^{n,i},\xi_h^{n,i}) \in \bm{V}_h \times W_h$, such that
  \begin{align}
    & a_1(\bm{u}_h^{n,i},\bm{v}_h)-b(\bm{v}_h,\xi_h^{n,i}) = ( \bm{f}^{n},\bm{v}_h ) + \langle \bm{f}_1^{n} , \bm{v}_h \rangle_{\Gamma_{\sigma}} , \quad \ \ \forall \bm{v}_h \in \bm{V}_h, \label{ts1}
    \\
    & b(\bm{u}_h^{n,i},\phi_h) + a_2(\xi_h^{n,i},\phi_h) = c(p_h^{n,i},\phi_h)  , \quad \quad \quad \quad \quad \quad \ \ \forall \phi_h \in W_h. \label{ts2}
  \end{align}
  \\
  \quad \quad \textbf{while} not converged.
  \\
  \quad \quad set $(\bm{u}_h^{n},\xi_h^{n},p_h^{n}) = (\bm{u}_h^{n,i},\xi_h^{n,i},p_h^{n,i})$.
  \\
  \textbf{end for}
\end{algorithm}

\begin{theorem}
\label{thm2}
For a fixed time step $n$, assume that the solution $(\bm{u}_h^{n-1}, \xi_h^{n-1}, p_h^{n-1})$, generated by \textbf{Algorithm \ref{algo:Coupled_CN}}, is provided as input to the equation \eqref{ts3} in \textbf{Algorithm \ref{algo:TIDA_CN}} to generate the sequence $\{(\bm{u}_h^{n,i},\xi_h^{n,i},p_h^{n,i})\}_{i \geq 0}$. 
The sequence converges to the solution $(\bm{u}_h^{n},\xi_h^{n},p_h^{n})$ produced by \textbf{Algorithm \ref{algo:Coupled_CN}}. The following estimate hold:
\begin{align}
       & \| e_{\xi,s}^{n,i} \|_{L^2(\Omega)} \leq L_s \| e_{\xi,s}^{n,i-1} \|_{L^2(\Omega)}, \label{thm2con} 
\end{align}
where $L_s$ is a positive constant strictly smaller than $1$ given by
\begin{align*}
    L_s = \frac{1}{ \left( \frac{c_0 \lambda}{\alpha^2} + \frac{k_p \Delta t \lambda}{2 C_P^2 \alpha^2} + 1  \right) \times \left( C(\Tilde{\beta},\mu) \lambda + 1 \right) }.
\end{align*}
Moreover, it holds that
\begin{align}
        & \| e_{p,s}^{n,i} \|_{H^1(\Omega)} \lesssim \| e_{\xi,s}^{n,i-1} \|_{L^2(\Omega)}, \label{thm2con2}\\ 
        & \|e_{\bm{u},s}^{n,i}\|_{H^1(\Omega)} \lesssim \|e_{\xi,s}^{n,i}\|_{L^2(\Omega)}.   \label{thm2con3}
    \end{align}
\end{theorem}

\begin{proof}
Subtracting \eqref{d1}, \eqref{d2}, \eqref{d3} from \eqref{ts1}, \eqref{ts2}, \eqref{ts3}, respectively, we see that
\begin{align}
    a_1(e_{\bm{u},s}^{n,i},\bm{v}_h)-b(\bm{v}_h,e_{\xi,s}^{n,i}) & = 0, \quad \quad \quad \quad \quad \quad \ \ \forall \bm{v}_h \in \bm{V}_h, \label{thm_ts:eq1}
    \\
    b(e_{\bm{u},s}^{n,i},\phi_h) + a_2(e_{\xi,s}^{n,i},\phi_h) & = c(e_{p,s}^{n,i},\phi_h), \quad \quad \quad \forall \phi_h \in W_h, \label{thm_ts:eq2}
    \\
    a_3 ( e_{p,s}^{n,i} , \psi_h  ) + \frac{\Delta t}{2}d ( e_{p,s}^{n,i},\psi_h ) & = c( \psi_h ,  e_{\xi,s}^{n,i-1} ), \quad \quad  \forall \psi_h \in M_h. \label{thm_ts:eq3}
\end{align}
Using the inf-sup condition \eqref{dinfsupcon}, we can derive the following inequality from \eqref{thm_ts:eq1}.
\begin{align}
    \Tilde{\beta} \| e_{\xi,s}^{n,i} \|_{L^2(\Omega)} & \leq \sup_{\bm{v}_h \in \bm{V}_h} \frac{  b(\bm{v}_h,e_{\xi,s}^{n,i}) }{\| \bm{v}_h \|_{H^1(\Omega)}} = \sup_{\bm{v}_h \in \bm{V}_h} \frac{  a_1(e_{\bm{u},s}^{n,i},\bm{v}_h) }{\| \bm{v}_h \|_{H^1(\Omega)}}
    \leq 2\mu C \| \varepsilon(e_{\bm{u},s}^{n,i}) \|_{L^2(\Omega)}, \label{thm_ts:eq4}
\end{align}
which leads to the following inequality
\begin{align}
    C(\Tilde{\beta},\mu) \| e_{\xi,s}^{n,i} \|_{L^2(\Omega)}^2 \leq 2\mu \| \varepsilon(e_{\bm{u},s}^{n,i}) \|_{L^2(\Omega)}^2. \label{thm_ts:eq4s}
\end{align}

For the reaction-diffusion problem, taking $\psi_h = e_{p,s}^{n,i}$ in \eqref{thm_ts:eq3} yields
\begin{align}
    a_3 ( e_{p,s}^{n,i} , e_{p,s}^{n,i}  ) + \frac{\Delta t}{2}d ( e_{p,s}^{n,i},e_{p,s}^{n,i} ) = c( e_{p,s}^{n,i} ,  e_{\xi,s}^{n,i-1} ). \label{thm_ts:eq5}
\end{align}
Applying the Poincar{\'e} inequality, we get
\begin{align}
    & \left( c_0 + \frac{\alpha^2}{\lambda} + \frac{k_p \Delta t}{2C_P^2} \right) \| e_{p,s}^{n,i} \|_{L^2(\Omega)}^2 \leq a_3 ( e_{p,s}^{n,i} , e_{p,s}^{n,i}  ) + \frac{\Delta t}{2}d ( e_{p,s}^{n,i},e_{p,s}^{n,i} ). \label{thm_ts:eq5s}
\end{align}
Using Cauchy-Schwarz inequality, there holds
\begin{align}
    & c( e_{p,s}^{n,i} ,  e_{\xi,s}^{n,i-1} ) \leq \frac{\alpha}{\lambda} \| e_{p,s}^{n,i} \|_{L^2(\Omega)} \| e_{\xi,s}^{n,i-1} \|_{L^2(\Omega)}. \label{thm_ts:eq6s}
\end{align}
Combining the equations \eqref{thm_ts:eq5}, \eqref{thm_ts:eq5s} and \eqref{thm_ts:eq6s}, we can deduce that
\begin{align}
    \left( c_0 + \frac{\alpha^2}{\lambda} + \frac{k_p \Delta t}{2C_P^2} \right) \| e_{p,s}^{n,i} \|_{L^2(\Omega)} \leq \frac{\alpha}{\lambda}  \| e_{\xi,s}^{n,i-1} \|_{L^2(\Omega)}. \label{thm_ts:eq6}
\end{align}

For the generalized Stokes problem, we take $\bm{v}_h = e_{\bm{u},s}^{n,i}$ in \eqref{thm_ts:eq1}, $\phi_h = e_{\xi,s}^{n,i}$ in \eqref{thm_ts:eq2} to derive
\begin{align}
    a_1(e_{\bm{u},s}^{n,i}, e_{\bm{u},s}^{n,i})-b( e_{\bm{u},s}^{n,i},e_{\xi,s}^{n,i}) & = 0, \label{thm_ts:eq1s}
    \\
    b(e_{\bm{u},s}^{n,i},e_{\xi,s}^{n,i}) + a_2(e_{\xi,s}^{n,i},e_{\xi,s}^{n,i}) & = c(e_{p,s}^{n,i},e_{\xi,s}^{n,i}). \label{thm_ts:eq2s}
\end{align}
Summing up the resulted equations \eqref{thm_ts:eq1s} and \eqref{thm_ts:eq2s} yields
\begin{align}
    a_1(e_{\bm{u},s}^{n,i},e_{\bm{u},s}^{n,i}) + a_2(e_{\xi,s}^{n,i},e_{\xi,s}^{n,i}) & = c(e_{p,s}^{n,i},e_{\xi,s}^{n,i}). \label{thm_ts:eq7}
\end{align}
By using the Cauchy-Schwarz inequality, we obtain
\begin{align}
    c(e_{p,s}^{n,i},e_{\xi,s}^{n,i}) \leq \frac{\alpha}{\lambda} \| e_{p,s}^{n,i} \|_{L^2(\Omega)} \| e_{\xi,s}^{n,i} \|_{L^2(\Omega)}. \label{thm_ts:eq7s}
\end{align}
Alongside \eqref{thm_ts:eq4s}, we can utilize \eqref{thm_ts:eq7s} to rephrase equation \eqref{thm_ts:eq7} in the following manner. 
\begin{align}
     \left( C(\Tilde{\beta},\mu) + \frac{1}{\lambda} \right) \| e_{\xi,s}^{n,i} \|_{L^2(\Omega)} \leq \frac{\alpha}{\lambda} \| e_{p,s}^{n,i} \|_{L^2(\Omega)}. \label{thm_ts:eq8}
\end{align}
Then, we combine \eqref{thm_ts:eq6} and \eqref{thm_ts:eq8} to obtain our main conclusion as expressed in \eqref{thm2con}. 

On the other hand, from \eqref{thm_ts:eq5} and \eqref{thm_ts:eq6s} we can deduce that
\begin{align}
    \frac{k_p \Delta t}{2} \|\nabla e_{p,s}^{n,i}\|_{L^2(\Omega)}^2 = \frac{\Delta t}{2} d( e_{p,s}^{n,i}, e_{p,s}^{n,i} ) \leq \frac{\alpha}{\lambda} \| e_{p,s}^{n,i} \|_{L^2(\Omega)} \| e_{\xi,s}^{n,i-1} \|_{L^2(\Omega)}. \label{thm_ts:eq9}
\end{align}
Applying the Poincar{\'e} inequality to \eqref{thm_ts:eq9} yields the conclusion (\ref{thm2con2}). 
Besides, from \eqref{thm_ts:eq1s}, we can deduce that
\begin{align}
    2 \mu \| \varepsilon(e_{\bm{u},s}^{n,i})\|_{L^2(\Omega)}^2 =  a_1(e_{\bm{u},s}^{n,i}, e_{\bm{u},s}^{n,i}) = b( e_{\bm{u},s}^{n,i},e_{\xi,s}^{n,i}) \leq \|\text{div} e_{\bm{u},s}^{n,i}\|_{L^2(\Omega)} \|e_{\xi,s}^{n,i}\|_{L^2(\Omega)}. \label{thm_ts:eq10}
\end{align}
Applying \eqref{korncon} to \eqref{thm_ts:eq10} yields (\ref{thm2con3}). 
This completes the proof.
\end{proof}

\begin{remark}
\label{remark1}
The input $(\bm{u}_h^{n-1}, \xi_h^{n-1}, p_h^{n-1})$, generated by \textbf{Algorithm \ref{algo:Coupled_CN}}, ensures that the approximations $(\bm{u}_h^{n,i}, \xi_h^{n,i}, p_h^{n,i})$ produced by \textbf{Algorithm \ref{algo:TIDA_CN}} converge globally to the numerical solution $(\bm{u}_h^{n}, \xi_h^{n}, p_h^{n})$ obtained from \textbf{Algorithm \ref{algo:Coupled_CN}}. Combined with Theorem \ref{CoThm}, this shows that the proposed Crank-Nicolson iterative decoupled algorithm achieves second-order accuracy in time.
In practice, we may select a sufficiently large iteration number $I$ and treat $(\bm{u}_h^{n}, \xi_h^{n}, p_h^{n}) \approx (\bm{u}_h^{n,I}, \xi_h^{n,I}, p_h^{n,I})$. However, this approach may introduce challenges related to error accumulation over time.
\end{remark}

\subsection{A global-in-time iterative decoupled algorithm}
\label{section_gida}

\begin{figure}[h]
\centering
\includegraphics[width=0.95\textwidth,trim=0 580 0 0,clip]{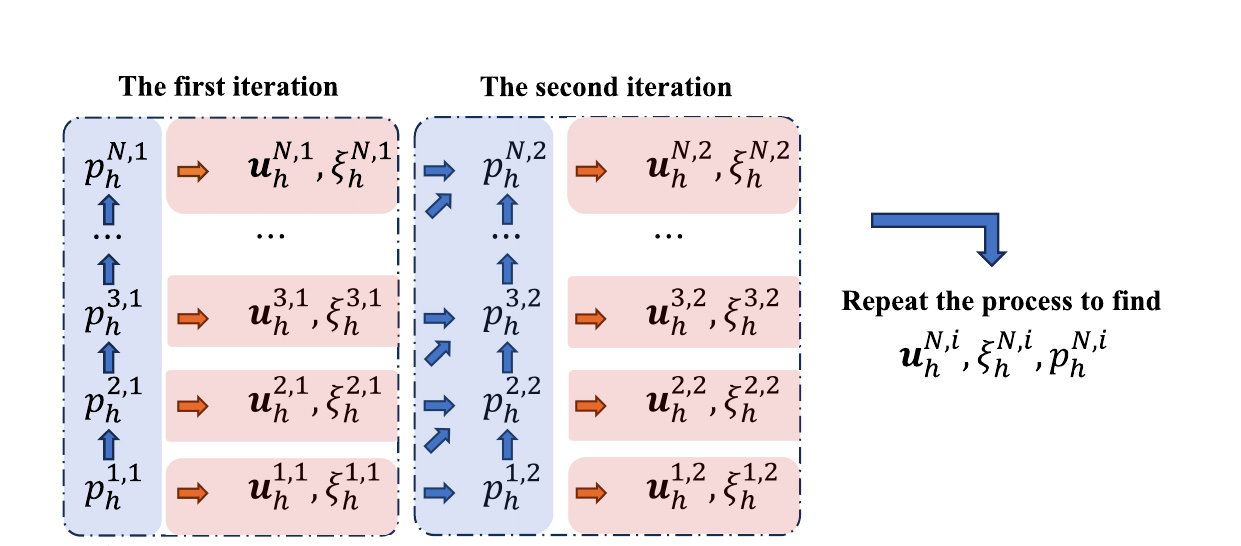}
\caption{
Progression in the global-in-time iterative decoupled algorithm.}
\label{fig_gt}
\end{figure}

The global-in-time method for poroelasticity, known for its compatibility with parallel processing, has seen recent attention \cite{borregales2019partially,ahmed2020adaptive}. Unlike traditional time-stepping, this approach solves subproblems across the entire time domain $[0, T]$.
We introduce a global-in-time Crank-Nicolson iterative decoupled algorithm in \textbf{Algorithm \ref{algo:GIDA_CN}}. For a visual representation, refer to the flowchart in Figure \ref{fig_gt}, where each iteration starts by sequentially solving the reaction-diffusion problem (blue zone), followed by the parallel solution of the partitioned generalized Stokes problem (red zone).
 
\begin{algorithm}[h]
\raggedright
  \caption{: A global-in-time Crank-Nicolson iterative decoupled algorithm}
  \label{algo:GIDA_CN}
  \textbf{Input:} initial information $(\bm{u}_h^{0},\xi_h^{0},p_h^{0}) \in \bm{V}_h \times W_h \times M_h$.
  \\
  \textbf{Output:} solution at the final time $(\bm{u}_h^{N},\xi_h^{N},p_h^{N}) \in \bm{V}_h \times W_h \times M_h$.
  \\
  set $i=0$, $(\bm{u}_h^{n,0},\xi_h^{n,0},p_h^{n,0})=(\bm{u}_h^{0},\xi_h^{0},p_h^{0})$, $n=0,1,\cdots, N$.
  \\
  \textbf{do} 
  \\
  \quad \quad  set $i = i + 1$,  $(\bm{u}_h^{0,i},\xi_h^{0,i},p_h^{0,i})=(\bm{u}_h^{0},\xi_h^{0},p_h^{0})$.
  \\
  \quad \quad \textbf{Step a:} \textbf{for} $n$ from $1$ to $N$
  \\
  \quad \quad \quad \quad find $p_h^{n,i}\in M_h$ such that
  \begin{align}
    & a_3 \left( \frac{p_h^{n,i} - p_h^{n-1,i}}{\Delta t}, \psi_h \right )-c\left( \psi_h , \frac{ \xi_h^{n,i-1}-\xi_h^{n-1,i-1} }{\Delta t} \right) + d \left( \frac{p_h^{n,i}+p_h^{n-1,i}}{2},\psi_h \right)
    \nonumber \\
    & \quad \quad = \frac{1}{2} ( Q_f^{n} + Q_f^{n-1} , \psi_h ) + \frac{1}{2} \langle  g_2^{n} + g_2^{n-1}  , \psi_h \rangle_{\Gamma_q}, \quad \quad \quad \forall \psi_h \in M_h. \label{gt3}
  \end{align}
  \\
  \quad \quad \textbf{end for}
  \\
  \quad \quad \textbf{Step b:} \textbf{parallel for} $n$ from $1$ to $N$
  \\
  \quad \quad \quad \quad find $(\bm{u}_h^{n,i},\xi_h^{n,i}) \in \bm{V}_h \times W_h$ such that
  \begin{align}
    & a_1(\bm{u}_h^{n,i},\bm{v}_h)-b(\bm{v}_h,\xi_h^{n,i}) = ( \bm{f}^{n},\bm{v}_h ) + \langle \bm{h}^{n} , \bm{v}_h \rangle_{\Gamma_{\sigma}} , \quad \ \ \forall \bm{v}_h \in \bm{V}_h, \label{gt1}
    \\
    & b(\bm{u}_h^{n,i},\phi_h) + a_2(\xi_h^{n,i},\phi_h) = c(p_h^{n,i},\phi_h)  , \quad \quad \quad \quad \quad \quad \ \ \forall \phi_h \in W_h. \label{gt2}
  \end{align}
  \\
  \quad \quad \textbf{end parallel for}
  \\
  \textbf{while} not converged.
  \\
  set $(\bm{u}_h^{N},\xi_h^{N},p_h^{N}) = (\bm{u}_h^{N,i},\xi_h^{N,i},p_h^{N,i})$.
\end{algorithm}

The convergence analysis of the global-in-time algorithm is provided in Theorem \ref{thm3}, establishing the relationship between \textbf{Algorithm \ref{algo:GIDA_CN}} and \textbf{Algorithm \ref{algo:Coupled_CN}}.
Let $(\bm{u}_h^n, \xi_h^n, p_h^n)$ and $(\bm{u}_h^{n,i}, \xi_h^{n,i}, p_h^{n,i})$ represent the solutions of problem \eqref{d1}-\eqref{d3} and problem \eqref{gt3}-\eqref{gt2}, respectively. We also define the $i$-th iteration errors as $e_{\bm{u},g}^{n,i} = \bm{u}_h^{n,i} - \bm{u}_h^n$, $e_{\xi,g}^{n,i} = \xi_h^{n,i} - \xi_h^n$, $e_{p,g}^{n,i} = p_h^{n,i} - p_h^n$.

\begin{theorem}
\label{thm3}
The sequence of solutions $\{(\bm{u}_h^{n,i}, \xi_h^{n,i}, p_h^{n,i})\}_{i \geq 0, 0 \leq n \leq N}$ generated by \textbf{ Algorithm \ref{algo:GIDA_CN}} converges to the solution $(\bm{u}_h^{n},\xi_h^{n},p_h^{n})$ produced by \textbf{Algorithm \ref{algo:Coupled_CN}}. The following estimate hold:
    \begin{align}
        \sum_{n=1}^N \| e_{\xi,g}^{n,i} - e_{\xi,g}^{n-1,i} \|_{L^2(\Omega)}^2 \leq L_g
        \sum_{n=1}^N \| e_{\xi,g}^{n,i-1} - e_{\xi,g}^{n-1,i-1} \|_{L^2(\Omega)}^2, \label{thm4con}
    \end{align}
where $L_g$ is a positive constant strictly smaller than 1 given by
\begin{align*}
    L_g = \frac{1}{ \left( \frac{c_0 \lambda}{\alpha^2} + 1  \right) \times \left( C(\beta,\mu) \lambda + 1 \right) }.
\end{align*}
Moreover, it holds that
    \begin{align}
        & \| e_{p,g}^{N,i} \|_{H^1(\Omega)}^2 \lesssim \sum_{n=1}^N \| e_{\xi,g}^{n,i-1} - e_{\xi,g}^{n-1,i-1} \|_{L^2(\Omega)}^2, \label{thm4con1} \\
        & \|e_{\xi,g}^{N,i}\|_{L^2(\Omega)} \lesssim \|e_{\bm{u},g}^{N,i}\|_{H^1(\Omega)} \lesssim   \|e_{p,g}^{N,i}\|_{H^1(\Omega)}.
    \end{align}
\end{theorem}

\begin{proof}
We begin by subtracting \eqref{d1}, \eqref{d2}, \eqref{d3} from \eqref{gt1}, \eqref{gt2}, \eqref{gt3}, respectively.
\begin{align}
    a_1(e_{\bm{u},g}^{n,i},\bm{v}_h)-b(\bm{v}_h,e_{\xi,g}^{n,i}) & = 0, \quad \quad \quad \quad \quad \quad \ \ \forall \bm{v}_h \in \bm{V}_h, \label{thm_gt:eq1}
    \\
    b(e_{\bm{u},g}^{n,i},\phi_h) + a_2(e_{\xi,g}^{n,i},\phi_h) & = c(e_{p,g}^{n,i},\phi_h), \quad \quad \quad \forall \phi_h \in W_h, \label{thm_gt:eq2}
    \\
    a_3 ( e_{p,g}^{n,i} - e_{p,g}^{n-1,i} , \psi_h  ) + \frac{\Delta t}{2}d ( e_{p,g}^{n,i} + e_{p,g}^{n-1,i},\psi_h ) & = c( \psi_h ,  e_{\xi,g}^{n,i-1} - e_{\xi,g}^{n-1,i-1} ), \quad  \forall \psi_h \in M_h. \label{thm_gt:eq3}
\end{align}

For the reaction-diffusion problem, we take $\psi_h = e_{p,g}^{n,i} - e_{p,g}^{n-1,i}$ in \eqref{thm_gt:eq3} to obtain
\begin{align}
     & \left( c_0 + \frac{\alpha^2}{\lambda} \right) \| e_{p,g}^{n,i} - e_{p,g}^{n-1,i} \|_{L^2(\Omega)}^2 + \frac{k_p \Delta t}{2} \| \nabla e_{p,g}^{n,i} \|_{L^2(\Omega)}^2 - \frac{k_p \Delta t}{2} \| \nabla e_{p,g}^{n-1,i} \|_{L^2(\Omega)}^2 
     \nonumber \\
     & = c( e_{p,g}^{n,i} - e_{p,g}^{n-1,i} ,  e_{\xi,g}^{n,i-1} - e_{\xi,g}^{n-1,i-1} ).  \label{thm_gt:eq7}
\end{align}
After using the Cauchy-Schwarz inequality on the right-hand side of \eqref{thm_gt:eq7}, we then apply the operator $\sum_{n=1}^N$ to the result to derive
\begin{align*}
    & \left( c_0 + \frac{\alpha^2}{\lambda} \right) \sum_{n=1}^N \| e_{p,g}^{n,i} - e_{p,g}^{n-1,i} \|_{L^2(\Omega)}^2 + \frac{k_p \Delta t}{2} \| \nabla e_{p,g}^{N,i} \|_{L^2(\Omega)}^2 
    \nonumber \\
    & \leq \frac{1}{2} \left( c_0 + \frac{\alpha^2}{\lambda} \right) \sum_{n=1}^N \| e_{p,g}^{n,i} - e_{p,g}^{n-1,i} \|_{L^2(\Omega)}^2 + \frac{1}{2 \left( c_0 + \frac{\alpha^2}{\lambda} \right) } \frac{\alpha^2}{\lambda^2} \sum_{n=1}^N \| e_{\xi,g}^{n,i-1} - e_{\xi,g}^{n-1,i-1} \|_{L^2(\Omega)}^2.
\end{align*}
The above inequality holds since the term $e_{p,g}^{0,i} = 0$, which leads directly to
\begin{align}
    & \sum_{n=1}^N \| e_{p,g}^{n,i} - e_{p,g}^{n-1,i} \|_{L^2(\Omega)}^2 + \frac{k_p \Delta t}{\left( c_0 + \frac{\alpha^2}{\lambda} \right)} \| \nabla e_{p,g}^{N,i} \|_{L^2(\Omega)}^2 
    \nonumber \\
    & \leq \frac{\frac{\alpha^2}{\lambda^2}}{\left( c_0 + \frac{\alpha^2}{\lambda} \right)^2 } \sum_{n=1}^N \| e_{\xi,g}^{n,i-1} - e_{\xi,g}^{n-1,i-1} \|_{L^2(\Omega)}^2.
    \label{thm_gt:eq8}
\end{align}

For the generalized Stokes problem, we take the difference of the $(n - 1)$-st time step and the $n$-th time step of \eqref{thm_gt:eq1} and \eqref{thm_gt:eq2}, respectively.
\begin{align}
    a_1(e_{\bm{u},g}^{n,i} - e_{\bm{u},g}^{n-1,i},\bm{v}_h)-b(\bm{v}_h,e_{\xi,g}^{n,i} - e_{\xi,g}^{n-1,i}) & = 0, \quad \quad \quad \quad \quad \quad \quad \quad \forall \bm{v}_h \in \bm{V}_h, \label{thm_gt:eq4}
    \\
    b(e_{\bm{u},g}^{n,i} - e_{\bm{u},g}^{n-1,i},\phi_h) + a_2(e_{\xi,g}^{n,i} - e_{\xi,g}^{n-1,i},\phi_h) & = c(e_{p,g}^{n,i} - e_{p,g}^{n-1,i},\phi_h), \quad \forall \phi_h \in W_h. \label{thm_gt:eq5}
\end{align}
Choosing $\bm{v}_h = e_{\bm{u},g}^{n,i} - e_{\bm{u},g}^{n-1,i}$ in \eqref{thm_gt:eq4}, $\phi_h = e_{\xi,g}^{n,i} - e_{\xi,g}^{n-1,i}$ in \eqref{thm_gt:eq5}, we then sum up the resulted equations to derive
\begin{align}
    2\mu \| \varepsilon (e_{\bm{u},g}^{n,i} - e_{\bm{u},g}^{n-1,i}) \|_{L^2(\Omega)}^2 + \frac{1}{\lambda} \| e_{\xi,g}^{n,i} - e_{\xi,g}^{n-1,i} \|_{L^2(\Omega)}^2 & = c(e_{p,g}^{n,i} - e_{p,g}^{n-1,i},e_{\xi,g}^{n,i} - e_{\xi,g}^{n-1,i}).  \label{thm_gt:eq6}
\end{align}
Using the inf-sup condition \eqref{dinfsupcon}, we can derive the following inequality from \eqref{thm_gt:eq4}.
\begin{align*}
    \Tilde{\beta} \| e_{\xi,g}^{n,i} - e_{\xi,g}^{n-1,i} \|_{L^2(\Omega)} & \leq \sup_{\bm{v}_h \in \bm{V}_h} \frac{  b(\bm{v}_h,e_{\xi,g}^{n,i} - e_{\xi,g}^{n-1,i} ) }{\| \bm{v}_h \|_{H^1(\Omega)}} \nonumber \\
    & = \sup_{\bm{v}_h \in \bm{V}_h} \frac{  a_1(e_{\bm{u},g}^{n,i} - e_{\bm{u},g}^{n-1,i},\bm{v}_h) }{\| \bm{v}_h \|_{H^1(\Omega)}}
    \leq 2\mu C \| \varepsilon(e_{\bm{u},g}^{n,i} - e_{\bm{u},g}^{n-1,i}) \|_{L^2(\Omega)}, 
\end{align*}
leading to the following inequality
\begin{align}
    C(\Tilde{\beta},\mu) \| e_{\xi,g}^{n,i} - e_{\xi,g}^{n-1,i} \|_{L^2(\Omega)}^2 \leq 2\mu \| \varepsilon(e_{\bm{u},g}^{n,i} - e_{\bm{u},g}^{n-1,i}) \|_{L^2(\Omega)}^2.
    \label{thm_gt:inf}
\end{align}
After applying the Cauchy-Schwarz inequality to the right-hand side of \eqref{thm_gt:eq6}, we then use \eqref{thm_gt:inf} to obtain that
\begin{align*}
     \left( C(\Tilde{\beta},\mu) + \frac{1}{\lambda} \right) \| e_{\xi,g}^{n,i} - e_{\xi,g}^{n-1,i} \|_{L^2(\Omega)} \leq \frac{\alpha}{\lambda} \| e_{p,g}^{n,i} - e_{p,g}^{n-1,i} \|_{L^2(\Omega)}.
\end{align*}
After taking the square of the above inequality, we then apply the summation operator $\sum_{n=1}^N$ to the result to obtain
\begin{align}
     \left( C(\Tilde{\beta},\mu) + \frac{1}{\lambda} \right)^2 \sum_{n=1}^N \| e_{\xi,g}^{n,i} - e_{\xi,g}^{n-1,i} \|_{L^2(\Omega)}^2 \leq \frac{\alpha^2} {\lambda^2} \sum_{n=1}^N \| e_{p,g}^{n,i} - e_{p,g}^{n-1,i} \|_{L^2(\Omega)}^2. 
     \label{thm_gt:eq6s}
\end{align}

Combining \eqref{thm_gt:eq8} and \eqref{thm_gt:eq6s} yields
\begin{align}
     & \frac{\lambda^2}{\alpha^2} \left( C(\Tilde{\beta},\mu) + \frac{1}{\lambda} \right)^2 \sum_{n=1}^N \| e_{\xi,g}^{n,i} - e_{\xi,g}^{n-1,i} \|_{L^2(\Omega)}^2 + \frac{k_p \Delta t}{\left( c_0 + \frac{\alpha^2}{\lambda} \right)} \| \nabla e_{p,g}^{N,i} \|_{L^2(\Omega)}^2 
    \nonumber \\
    & \leq \frac{\frac{\alpha^2}{\lambda^2}}{\left( c_0 + \frac{\alpha^2}{\lambda} \right)^2 } \sum_{n=1}^N \| e_{\xi,g}^{n,i-1} - e_{\xi,g}^{n-1,i-1} \|_{L^2(\Omega)}^2. \label{thm_gt:eqmc}
\end{align}
Here, we can easily deduce our main conclusions \eqref{thm4con} and \eqref{thm4con1} from the inequality \eqref{thm_gt:eqmc}. 
On the other hand, we can take $\bm{v}_h = e_{\bm{u},g}^{N,i}$ in the $N$-th step of \eqref{thm_gt:eq1}, and  $\phi_h = e_{\xi,g}^{N,i}$ in the $N$-th step of \eqref{thm_gt:eq2} to obtain
\begin{align}
    a_1(e_{\bm{u},g}^{N,i}, e_{\bm{u},g}^{N,i}) - b( e_{\bm{u},g}^{N,i},e_{\xi,g}^{N,i}) & = 0, \label{thm_gt:eq13} \\
    b(e_{\bm{u},g}^{N,i},e_{\xi,g}^{N,i}) + a_2(e_{\xi,g}^{N,i},e_{\xi,g}^{N,i}) & = c(e_{p,g}^{N,i},e_{\xi,g}^{N,i}). \label{thm_gt:eq14}    
\end{align}
By applying the Cauchy-Schwarz inequality to \eqref{thm_gt:eq13}, there holds
\begin{align}
    2\mu \| \varepsilon(e_{\bm{u},g}^{N,i}) \|_{L^2(\Omega)}^2 \leq \|\mbox{div}e_{\bm{u},g}^{N,i} \|_{L^2(\Omega)} \| e_{\xi,g}^{N,i} \|_{L^2(\Omega)}.
\end{align}
The above inequality implies $\| e_{\xi,g}^{N,i} \|_{L^2(\Omega)} \lesssim \|e_{\bm{u},g}^{N,i}\|_{H^1(\Omega)}$.
Summing up \eqref{thm_gt:eq13} and \eqref{thm_gt:eq14}, we can apply the Cauchy-Schwarz inequality to obtain
\begin{align*}
    \frac{1}{\lambda}  \| e_{\xi,g}^{N,i} \|_{L^2(\Omega)}^2 & \leq a_1(e_{\bm{u},g}^{N,i}, e_{\bm{u},g}^{N,i}) + a_2(e_{\xi,g}^{N,i},e_{\xi,g}^{N,i}) = c(e_{p,g}^{N,i},e_{\xi,g}^{N,i}) \leq \frac{\alpha}{\lambda} \|e_{p,g}^{N,i}\|_{L^2(\Omega)} \|e_{\xi,g}^{N,i}\|_{L^2(\Omega)},
\end{align*}
which leads to $\|e_{\bm{u},g}^{N,i}\|_{H^1(\Omega)} \lesssim \|e_{p,g}^{N,i}\|_{H^1(\Omega)}$. The proof is complete.
\end{proof}

From the above theorem, we observe that the error term  $\sum_{n=1}^N \| e_{\xi,g}^{n, i} - e_{\xi,g}^{n-1, i}\|_{L^2(\Omega)}^2$ consistently diminishes to zero, leading to the convergence of the other error terms, $\|e_{\bm{u},g}^{N, i}\|_{H^1(\Omega)}$, $\|e_{\xi,g}^{N,i}\|_{L^2(\Omega)}$, and $\|e_{p,g}^{N, i}\|_{H^1(\Omega)}$, also approaching zero. In alignment with Theorem \ref{CoThm}, it is evident that \textbf{Algorithm \ref{algo:GIDA_CN}} achieves second-order accuracy in time. 
However, it is important to note that the solution errors $\|e_{\bm{u},g}^{N, i}\|_{H^1(\Omega)}$, $\|e_{\xi,g}^{N,i}\|_{L^2(\Omega)}$, and $\|e_{p,g}^{N, i}\|_{H^1(\Omega)}$ may not decrease monotonically. In long-time simulations, where $T$ is large and a significant number of time steps are required, the global-in-time approach efficiently leverages parallel computing strategies. Further insights on parallel efficiency are provided below.

\begin{prop}[Speedup of the global-in-time algorithm]
\label{speedup} 
Assume that there are a total of $2^k = N$ time steps and $2^m$ processors. Let $T_p$ represent the CPU time required to compute $\{p_h^{n,i}\}_{1 \leq n \leq N}$ in the $i$-th iteration on a serial processor, and let $T_{(\bm{u},\xi)}$ denote the CPU time spent computing $\{(\bm{u}_h^{n,i}, \xi_h^{n,i})\}_{1 \leq n \leq N}$ in the $i$-th iteration on a serial processor. According to Amdahl’s law \cite{gustafson1988reevaluating}, the speedup is given by:
\begin{align*}
    Speedup = \frac{T_{p}+T_{(\bm{u},\xi)}}{T_{p} + 2^{-m} \times T_{(\bm{u},\xi)}},
\end{align*}
where $m \leq k$ is considered. 
\end{prop}

As the number of processors increases, the speedup becomes more pronounced, which is particularly true in cases where $T_{(\bm{u},\xi)}$ dominates the computation time.
This occurs because the solutions $\{\bm{u}_h^{n, i}\}$  represent some vector fields and thus require more computations, whereas the solutions $\{p_h^{n, i}\}$ are some scalar fields and involve fewer operations. This distinction renders the global-in-time algorithm highly suitable for parallel implementation, thereby enabling significant efficiency gains in large-scale simulations. It is particularly advantageous in cases that entail long-time integration and fine spatial-temporal resolution, where computational demands are high.

\section{Numerical experiments}
\label{Sec5}
In this section, we present three two-dimensional numerical experiments to evaluate the convergence behavior of the proposed algorithms and validate our theoretic predictions. In the first example, we examine the second-order time accuracy of \textbf{Algorithm \ref{algo:Coupled_CN}} based on a model with an analytic solution. The second example is the Barry-Mercer model, which has a singularity in the source term and the physical parameters include some degenerate cases (for example, $c_0=0$ and the permeability coefficient is small, $K = 10^{-6}$).  The well-known benchmark Mandel's problem is explored in the third example.  The computations are performed using the open-source finite element software FEniCS \cite{alnaes2015fenics}. 

For the iterative decoupled algorithms, let $I$ be the maximum number of iterations. At the final time $T = t_N$, the exact solution is denoted as $(\bm{u}^{N}, \xi^{N}, p^{N})$. We denote the output of \textbf{Algorithm \ref{algo:Coupled_CN}}, \textbf{Algorithm \ref{algo:TIDA_CN}}, \textbf{Algorithm \ref{algo:GIDA_CN}} as $(\bm{u}_{h}^{N}, \xi_{h}^{N}, p_{h}^{N})$, $(\bm{u}_{h,s}^{N,I}, \xi_{h,s}^{N,I}, p_{h,s}^{N,I})$, $(\bm{u}_{h,g}^{N,I}, \xi_{h,g}^{N,I}, p_{h,g}^{N,I})$, respectively. The errors for \textbf{Algorithm \ref{algo:TIDA_CN}} can be decomposed into two parts, as outlined below:
\begin{align*}
    & \bm{u}_{h,s}^{N,I} - \bm{u}^N = (\bm{u}_{h,s}^{N,I} - \bm{u}_{h}^{N}) + ( \bm{u}_{h}^{N} - \bm{u}^N ) := e_{\bm{u},s}^{N,I} + ( \bm{u}_{h}^{N} - \bm{u}^N ),
    \\
    & \xi_{h,s}^{N,I} - \xi^N = (\xi_{h,s}^{N,I} - \xi_{h}^{N}) + ( \xi_{h}^{N} - \xi^N ) := e_{\xi,s}^{N,I} + ( \xi_{h}^{N} - \xi^N ),
    \\
    & p_{h,s}^{N,I} - p^N = (p_{h,s}^{N,I} - p_{h}^{N}) + ( p_{h}^{N} - p^N ) := e_{p,s}^{N,I} + ( p_{h}^{N} - p^N ).
\end{align*}
Similarly, the errors for \textbf{Algorithm \ref{algo:GIDA_CN}} are expressed as:
\begin{align*}
    & \bm{u}_{h,g}^{N,I} - \bm{u}^N = (\bm{u}_{h,g}^{N,I} - \bm{u}_{h}^{N}) + ( \bm{u}_{h}^{N} - \bm{u}^N ) := e_{\bm{u},g}^{N,I} + ( \bm{u}_{h}^{N} - \bm{u}^N ),
    \\
    & \xi_{h,g}^{N,I} - \xi^N = (\xi_{h,g}^{N,I} - \xi_{h}^{N}) + ( \xi_{h}^{N} - \xi^N ) := e_{\xi,g}^{N,I} + ( \xi_{h}^{N} - \xi^N ),
    \\
    & p_{h,g}^{N,I} - p^N = (p_{h,g}^{N,I} - p_{h}^{N}) + ( p_{h}^{N} - p^N ) := e_{p,g}^{N,I} + ( p_{h}^{N} - p^N ).
\end{align*}
By Theorem \ref{CoThm}, the error terms $( \bm{u}_{h}^{N} - \bm{u}^N )$, $( \xi_{h}^{N} - \xi^N )$, and $( p_{h}^{N} - p^N )$ should exhibit second-order accuracy in time. By Theorem \ref{thm2} and Theorem \ref{thm3}, the error terms $(e_{\bm{u},s}^{N,I}, e_{\xi,s}^{N,I}, e_{p,s}^{N,I})$ and $(e_{\bm{u},g}^{N,I}, e_{\xi,g}^{N,I}, e_{p,g}^{N,I})$ should converge to $(\bm{0},0,0)$ when $I$ goes to infinity. We will verify all these theoretical predictions.

\subsection{Example 1: the convergence behavior of Algorithms \ref{algo:Coupled_CN}, \ref{algo:TIDA_CN}, \ref{algo:GIDA_CN}.}
\label{subsection4.1}
To confirm optimal convergence rates in time, we conduct the numerical experiment outlined in \cite{chaabane2018splitting}. 
The experiment considers a unit square domain $\Omega=[0,1]^2$, with a final time of $T=1.0$. 
Pure Dirichlet conditions are specified at $\partial \Omega$. 
The physical parameters are configured as follows:
\begin{align*}
    \lambda =1.0, \quad \mu = 1.0, \quad \alpha = 1.0, \quad c_0=1.0, \quad k_p=1.0.  
\end{align*}
The body force $\bm{f}$, source or sink term $Q_f$, initial conditions, and boundary conditions are chosen such that the exact solutions are given by:
\begin{align*}
    u_1 = \frac{1}{10} e^t (x+y^3), \quad
    u_2 = \frac{1}{10} t^2 (x^3+y^3), \quad
    p = 10e^{\frac{x+y}{10}}(1+t^3).
\end{align*} 
 
To assess the convergence rates in time, we utilize a small mesh size of $h = 1/256$ and refine the time step size $\Delta t$. The results, including errors and convergence rates obtained by \textbf{Algorithm \ref{algo:Coupled_CN}}, are presented in Table \ref{E1}. It is evident from the table that the error terms $( \bm{u}_{h}^{N} - \bm{u}^N )$, $( \xi_{h}^{N} - \xi^N )$, and $( p_{h}^{N} - p^N )$ display second-order convergence in time, consistent with the expectations from Theorem \ref{CoThm}.

\begin{table}[ht]
\setlength{\belowcaptionskip}{0.1cm}
\begin{center}
\renewcommand\arraystretch{1.5}
\caption{Errors and convergence rates of \textbf{Algorithm \ref{algo:Coupled_CN}} for Example 1.}
\label{E1}
\centering
{  \footnotesize
	\begin{tabular}{c|cl|cl|cl}
	\hline
	$\Delta t$ & \multicolumn{1}{l}{$\|\bm{u}_{h}^{N} - \bm{u}^N\|_{H^1(\Omega)}$} & Orders & \multicolumn{1}{l}{$\|\xi_{h}^{N} - \xi^N\|_{L^2(\Omega)}$} & Orders & \multicolumn{1}{l}{$\| p_{h}^{N} - p^N \|_{H^1(\Omega)}$} & \multicolumn{1}{l}{Orders} \\ \hline
    1/2  & 6.620e-03 &      & 3.070e-02 &      & 1.540e-01 &      \\
    1/4  & 2.629e-03 & 1.33 & 1.266e-02 & 1.28 & 6.333e-02 & 1.28 \\
    1/8  & 6.426e-04 & 2.03 & 3.297e-03 & 1.94 & 1.655e-02 & 1.94 \\
    \textbf{1/16} & \textbf{1.586e-04} & \textbf{2.02} & \textbf{8.285e-04} & \textbf{1.99} & \textbf{4.233e-03} & \textbf{1.97} \\
    \hline
	\end{tabular}
}
\end{center}
\end{table}

\begin{figure}
{\color{black}
\centering
\subfigure{\includegraphics[width=0.75\textwidth]{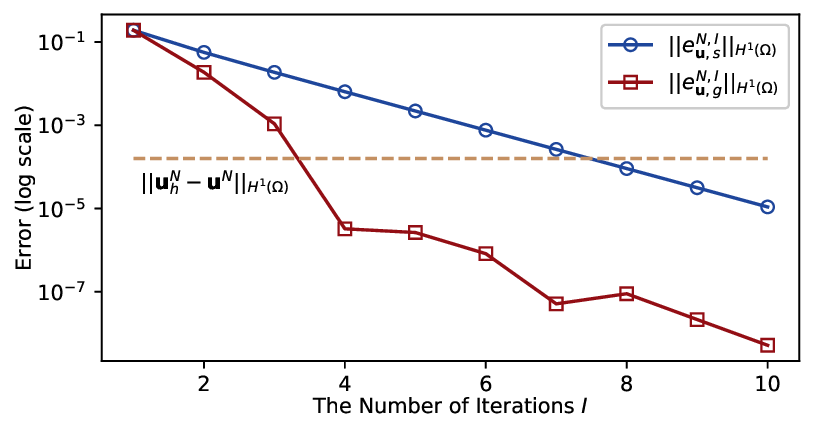}}
\subfigure{\includegraphics[width=0.75\textwidth]{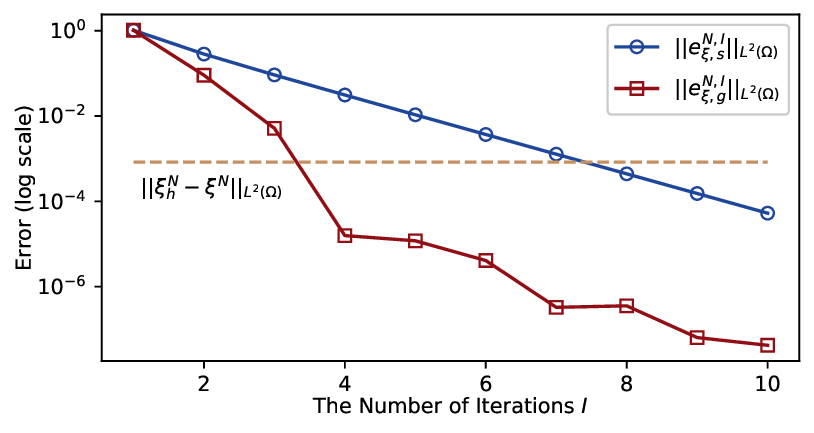}}
\subfigure{\includegraphics[width=0.75\textwidth]{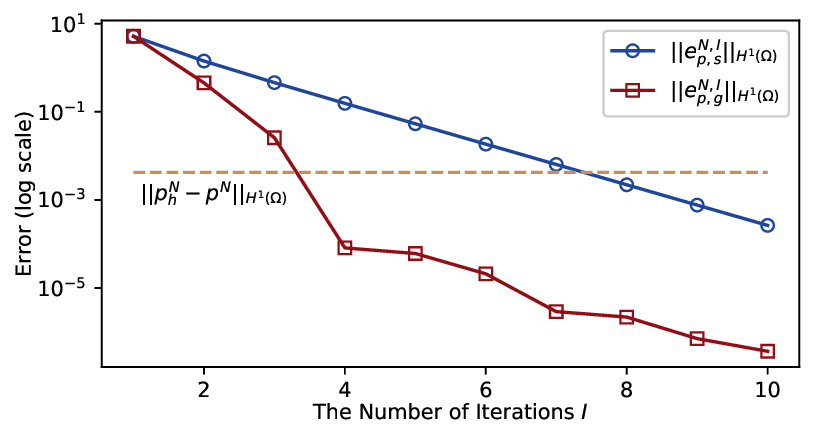}}
\caption{Convergence behaviors of \textbf{Algorithm \ref{algo:TIDA_CN}} and \textbf{Algorithm \ref{algo:GIDA_CN}} for Example 1.}
\label{E2}
}
\end{figure}

Using the same mesh and a fixed time step size of $\Delta t = 1/16$, we closely investigate the convergence behaviors of both \textbf{Algorithm \ref{algo:TIDA_CN}} and \textbf{Algorithm \ref{algo:GIDA_CN}}. Figure \ref{E2}  illustrates the convergence of the error terms $(e_{\bm{u},s}^{N,I}, e_{\xi,s}^{N,I}, e_{p,s}^{N,I})$ and $(e_{\bm{u,g}}^{N,I}, e_{\xi,g}^{N,I}, e_{p,g}^{N,I})$.
The results show that the error in both iterative decoupled algorithms decreases as the number of iterations $I$ increases. Notably, \textbf{Algorithm \ref{algo:TIDA_CN}} requires 8 iterations to ensure the error bounds:
\begin{align*}
&\| e_{\bm{u},s}^{N,8} \|_{H^1(\Omega)} \leq \|\bm{u}_{h}^{N} - \bm{u}^N\|_{H^1(\Omega)}, \\
&\| e_{\xi,s}^{N,8} \|_{L^2(\Omega)} \leq \|\xi_{h}^{N} - \xi^N\|_{L^2(\Omega)}, \\
&\| e_{p,s}^{N,8} \|_{H^1(\Omega)} \leq \|p_{h}^{N} - p^N\|_{H^1(\Omega)}.
\end{align*}
In comparison, \textbf{Algorithm \ref{algo:GIDA_CN}} achieves the same bounds with only 4 iterations:
\begin{align*}
&\| e_{\bm{u},g}^{N,4} \|_{H^1(\Omega)} \leq \|\bm{u}_{h}^{N} - \bm{u}^N\|_{H^1(\Omega)}, \\
&\| e_{\xi,g}^{N,4} \|_{L^2(\Omega)} \leq \|\xi_{h}^{N} - \xi^N\|_{L^2(\Omega)}, \\
&\| e_{p,g}^{N,4} \|_{H^1(\Omega)} \leq \|p_{h}^{N} - p^N\|_{H^1(\Omega)}.
\end{align*}
It is important to note that the oscillation observed in the error curve for \textbf{Algorithm \ref{algo:GIDA_CN}} is explainable. As discussed in Theorem \ref{thm3}, the term $\sum_{n=1}^N \| e_{\xi,g}^{n, i} - e_{\xi,g}^{n-1, i} \|_{L^2(\Omega)}^2$ only serves as an upper bound for the other error terms. While $\|e_{\bm{u},g}^{N, i}\|_{H^1(\Omega)}$ and $\|e_{p,g}^{N, i}\|_{H^1(\Omega)}$ converge, they may not decrease monotonically throughout the iterations.

\newpage
\subsection{Example 2: Barry-Mercer's model}
\label{subsection4.2}
The Barry-Mercer model, widely recognized in the literature \cite{phillips2005finite,barry1999exact}, serves as a benchmark with a point-source term. We consider the domain $\Omega=[0,1]^2$ where the point source is positioned at $(x_0, y_0) = (0.25, 0.25)$. The physical parameters are specified as follows \cite{cai2023some}.
\begin{align*}
    E=10^5, \quad \nu = 0.1, \quad \alpha = 1.0, \quad c_0=0.0, \quad k_p=10^{-6}, 
\end{align*}
which gives $\lambda = 1.1 \times 10^4$, $\mu = 4.5 \times 10^4$. The boundary segments are: 
$\Gamma_1=\{(1,y); 0\leq y\leq1\}$, $\Gamma_2=\{(x,0);0\leq x\leq1\}$, $\Gamma_3=\{(0,y);0\leq y\leq1\}$, $\Gamma_4=\{(x,1);0\leq x\leq1\}$. The initial and boundary conditions are specified as:
\begin{align*}
    \bm{u} = 0, \ p = 0 \quad & \mbox{in} \ \Omega \times \{0\}, \\
        \frac{\partial u_1}{\partial x} = 0, \ u_2 = 0 \quad & \mbox{on} \ \Gamma_j \times (0,T], \ j=1,3,  \\
	\frac{\partial u_2}{\partial y} = 0, \ u_1 = 0 \quad & \mbox{on} \ \Gamma_j \times (0,T], \ j=2,4,  \\
	p = 0, \ \bm{h} = \bm{0}, \ g_2 = 0 \quad & \mbox{on} \ \Gamma_j \times (0,T], \ j=1,2,3,4.
\end{align*}
The body force term $\bm{f} = \bm{0}$, and the source/sink term is given by:
\begin{align*}
    Q_f = 2 \omega \delta (x- x_0) \delta (y- y_0) \sin{(\omega t)},
\end{align*}
where $\omega = (\lambda + 2\mu) K$, and $\delta (\cdot)$ represents the Dirac function. 

Given that the analytic solutions for Barry-Mercer's problem are expressed as a series \cite{phillips2005finite}, obtaining exact values of $(\bm{u}^N, \xi^N, p^N)$ can be computationally challenging. To simplify, we focus on understanding the relationship among $(\bm{u}_{h}^{N}, \xi_{h}^{N}, p_{h}^{N})$, $(\bm{u}_{h,s}^{N,I}, \xi_{h,s}^{N,I}, p_{h,s}^{N,I})$ and $(\bm{u}_{h,g}^{N,I}, \xi_{h,g}^{N,I}, p_{h,g}^{N,I})$.
We employ the following relative errors as measures to check the convergence behaviors of \textbf{Algorithm \ref{algo:TIDA_CN}} and \textbf{Algorithm \ref{algo:GIDA_CN}}:
\begin{align*}
    & RE(\bm{u}_{h,s}^{N,I}) = \|e_{\bm{u},s}^{N,I}\|_{H^1(\Omega)} / \| \bm{u}_{h}^{N} \|_{H^1(\Omega)}, \quad RE(\bm{u}_{h,g}^{N,I}) = \|e_{\bm{u},g}^{N,I}\|_{H^1(\Omega)} / \| \bm{u}_{h}^{N} \|_{H^1(\Omega)}, 
    \\
    & RE(\xi_{h,s}^{N,I}) = \|e_{\xi,s}^{N,I}\|_{L^2(\Omega)} / \| \xi_{h}^{N} \|_{L^2(\Omega)}, \quad \ \ RE(\xi_{h,g}^{N,I}) = \|e_{\xi,g}^{N,I}\|_{L^2(\Omega)} / \| \xi_{h}^{N} \|_{L^2(\Omega)},
    \\
    & RE(p_{h,s}^{N,I}) = \|e_{p,s}^{N,I}\|_{H^1(\Omega)} / \| p_{h}^{N} \|_{H^1(\Omega)}, \quad \ RE(p_{h,g}^{N,I}) = \|e_{p,g}^{N,I}\|_{H^1(\Omega)} / \| p_{h}^{N} \|_{H^1(\Omega)}.
\end{align*}

\begin{figure}
\centering
\subfigure{\includegraphics[width=0.75\textwidth]{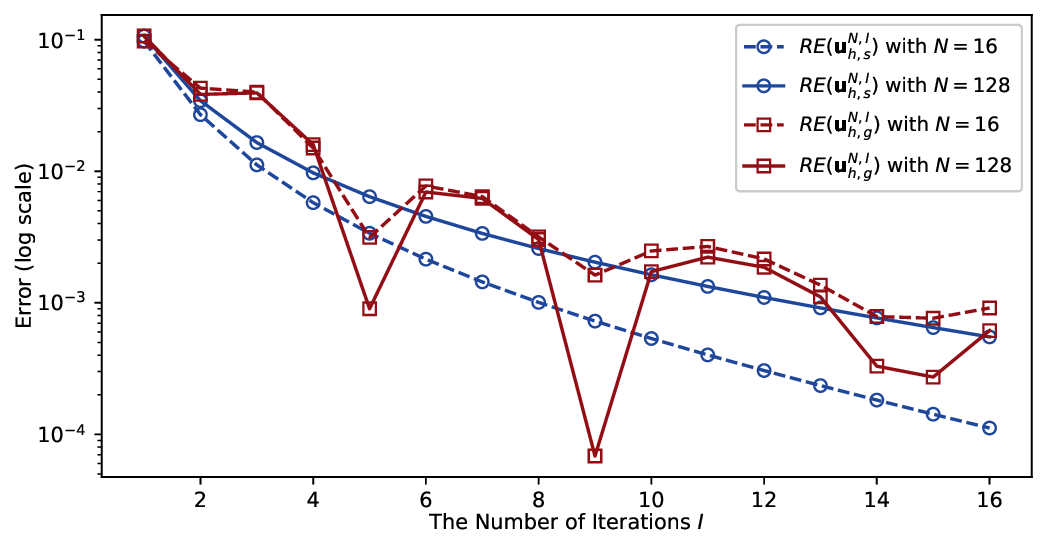}}
\subfigure{\includegraphics[width=0.75\textwidth]{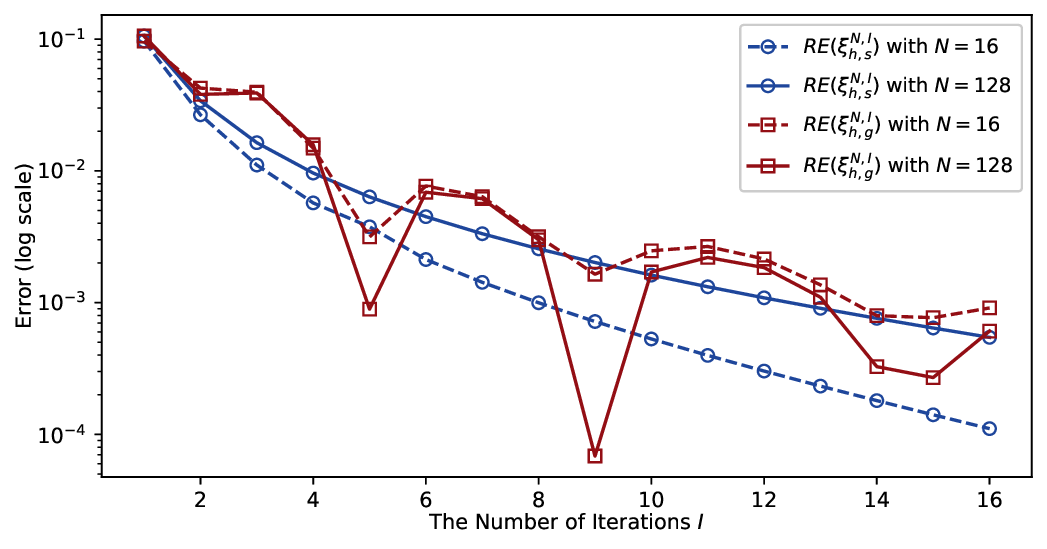}}
\subfigure{\includegraphics[width=0.75\textwidth]{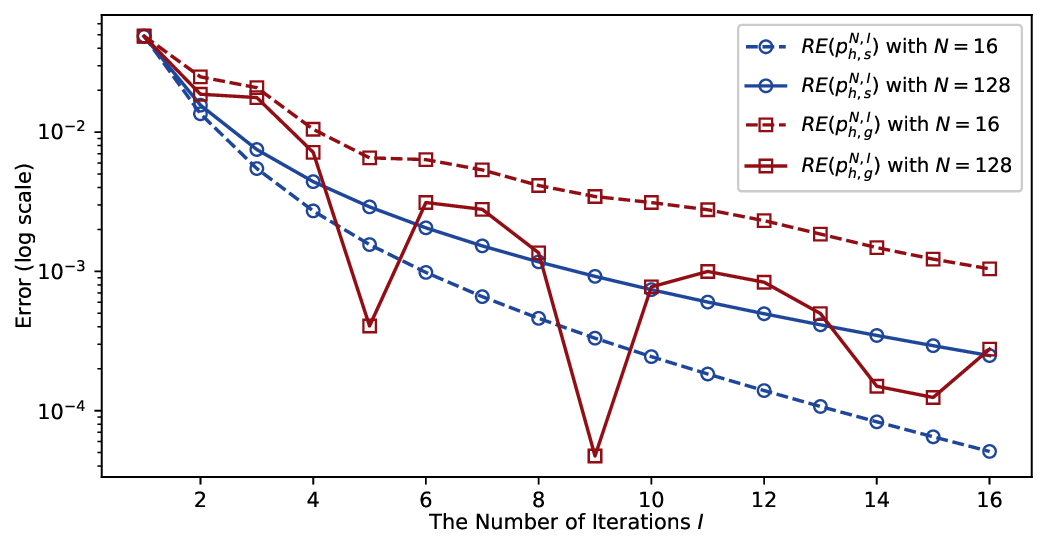}}
\caption{Convergence behaviors of \textbf{Algorithm \ref{algo:TIDA_CN}} and \textbf{Algorithm \ref{algo:GIDA_CN}} for Example 2.}
\label{E3}
\end{figure}

For our experiments, we adopt a mesh size of $h = 1/20$, keep the final time fixed at $T = \pi/(2\omega)$, and investigate two distinct time step sizes: $\Delta t = T/16$ and $\Delta t = T/128$. The numerical results are summarized in Figure \ref{E3}.
From Figure \ref{E3}, it is evident that the time-stepping iterative decoupled algorithm exhibits faster convergence than the global-in-time iterative decoupled algorithm for $N = 16$. However, for $N = 128$, the convergence rates of the two algorithms become comparable. Notably, as the time step size $\Delta t$ decreases, the convergence of the time-stepping iterative decoupled algorithm slows, indicating the presence of error accumulation. In contrast, the global-in-time iterative decoupled algorithm shows an opposing trend, with an accelerated convergence rate under smaller time step sizes. This observation suggests that while the time-stepping algorithm proves more efficient with larger time steps, the global-in-time algorithm emerges as a competitive alternative, particularly in scenarios involving smaller time step sizes.

\subsection{Example 3: Mandel's problem}

\begin{figure}
\centering
\subfigure{
\includegraphics[width=0.45\textwidth,trim=0 320 0 0,clip]{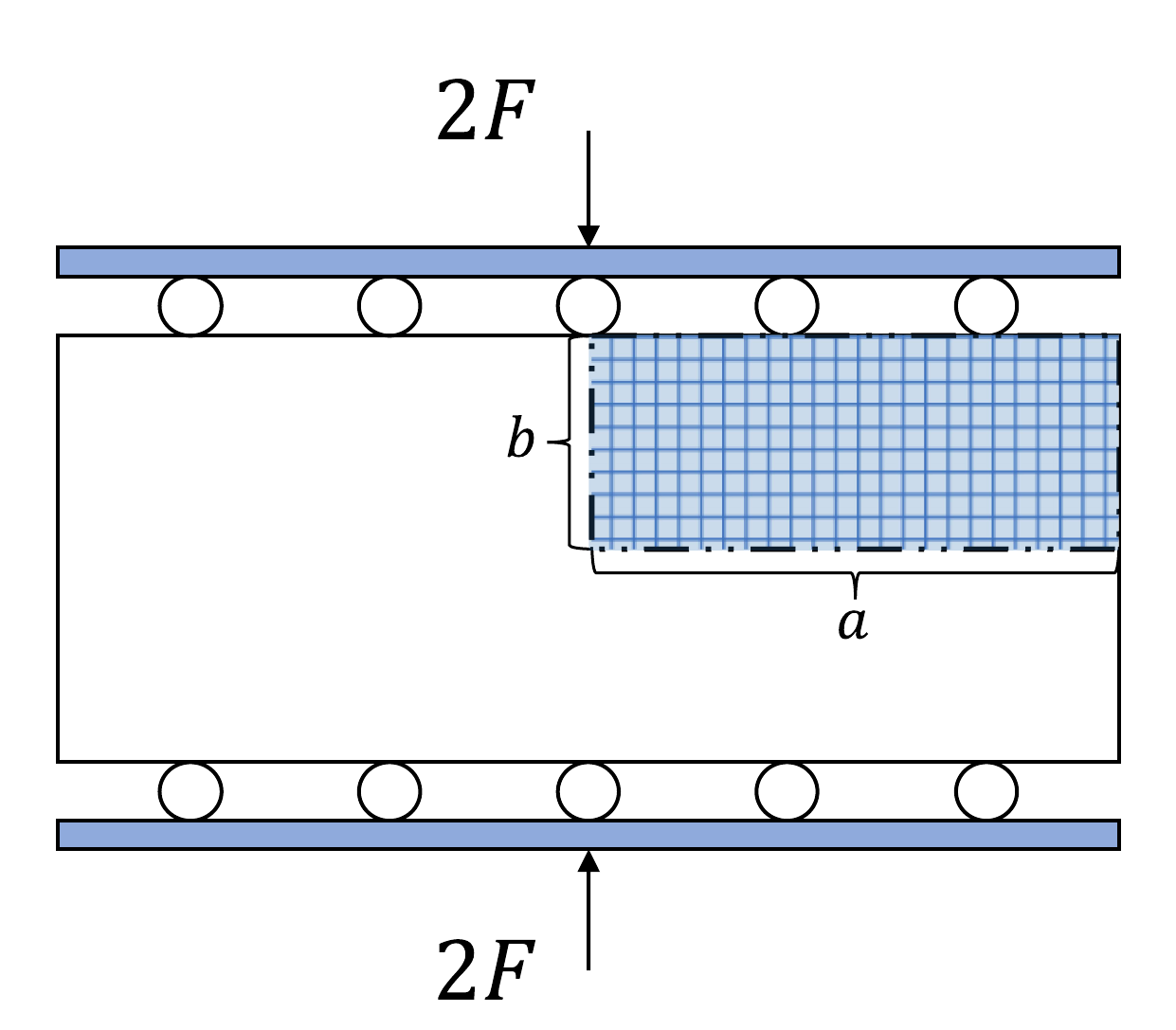}}
\subfigure{
\includegraphics[width=0.48\textwidth,trim= 40 420 20 20,clip]{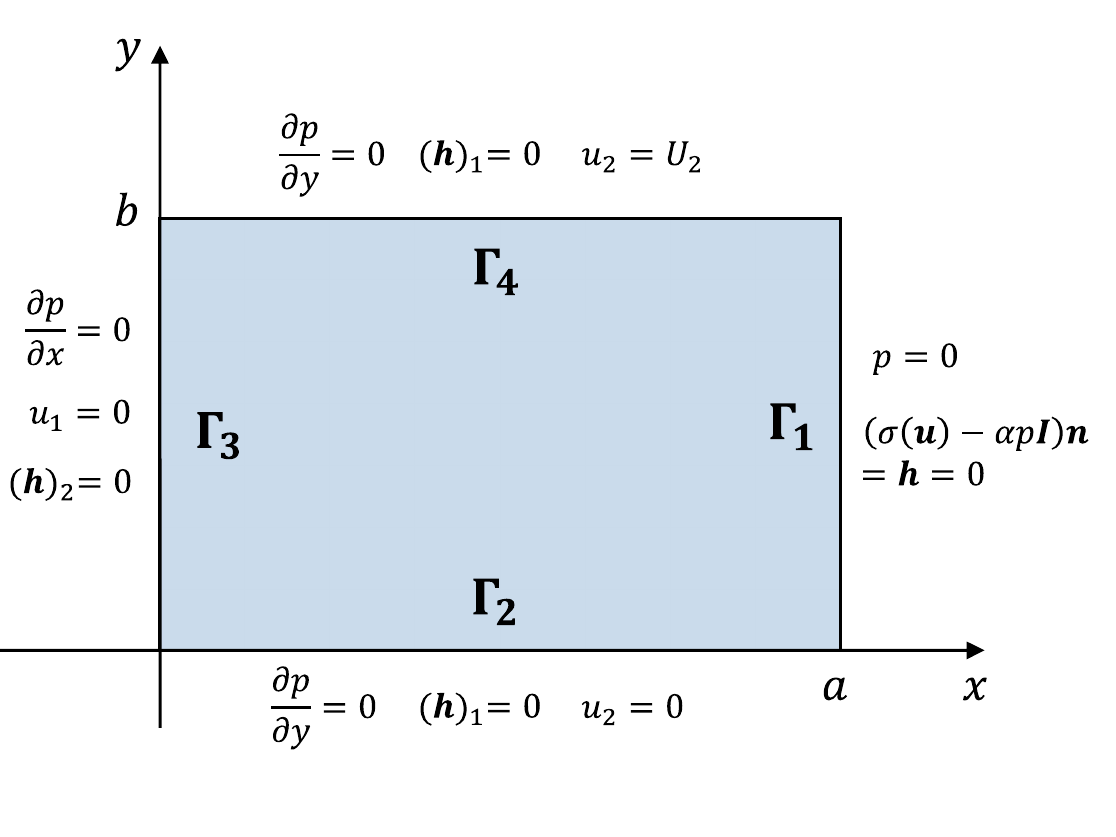}}
\caption{Description of Mandel's Problem: (Left) Computational domain within the blue zone and (Right) the corresponding boundary condition settings.}
\label{figMP}
\end{figure}

Mandel's problem \cite{mandel1953consolidation,abousleiman1996mandel} is a classical benchmark known for the Mandel-Cryer effect, and we use it to assess our proposed algorithms. This problem involves a saturated poroelastic plate in the domain $\Omega = [-a, a] \times [-b,b]$, constrained between two rigid plates (see the left side of Figure \ref{figMP}). The top plate exerts a constant downward force of $2F$, while the bottom plate applies an equal upward force. Due to the quarter symmetry of Mandel's problem, we simplify the computational domain to $[0, a] \times [0,b]$. For computational ease, we set $a = b = 1$, consistent with the boundary segments outlined in Section \ref{subsection4.2}. The right side of Figure \ref{figMP} illustrates the boundary conditions, following the approach used in \cite{mikelic2014numerical,almani2016convergence}.
\begin{align*}
        p = 0, \quad (\bm{h})_1 = (\bm{h})_2 = 0 \quad \ \ \quad & \mbox{on} \ \Gamma_1 \times (0,T],  \\
	\frac{\partial p}{\partial y} = 0, \quad (\bm{h})_1 =0, \quad u_2 = 0 \ \ \quad & \mbox{on} \ \Gamma_2 \times (0,T],  \\
        \frac{\partial p}{\partial x} = 0, \quad (\bm{h})_2 =0, \quad u_1 = 0 \ \ \quad & \mbox{on} \ \Gamma_3 \times (0,T],  \\
	\frac{\partial p}{\partial y} = 0, \quad (\bm{h})_1 =0, \quad u_2 = U_2 \quad & \mbox{on} \ \Gamma_4 \times (0,T],
\end{align*}
where $U_2$ is the closed-form solution for $y$-displacement. The body force and the source term are set to $\bm{f}=Q_f=0$. Following the experimental settings in \cite{phillips2005finite,phillips2007coupling,mikelic2014numerical}, we consider the following instantaneous pressure and deformation as the initial conditions:
\begin{align*}
    & u_1(0) = \lim_{t \rightarrow 0^+} U_1= \frac{F \nu_u}{2\mu}x, \\
    & u_2(0) = \lim_{t \rightarrow 0^+} U_2= \frac{-F(1-\nu_u)}{2\mu}y, \\ 
    & p(0)   = \lim_{t \rightarrow 0^+} P = \frac{FB(1+\nu_u)}{3} \quad \quad \mbox{in} \ \Omega \times \{0\},
\end{align*}
where $U_1$ is the closed-form solution for $x$-displacement, $P$ is the closed-form solution for pressure, $B$ is Skempton's coefficient, and $\nu_u$ is the undrained Poisson’s ratio. The detailed setting of the physical parameters is listed in Table \ref{PP}. 

\begin{table}
\setlength{\belowcaptionskip}{0.1cm}
\begin{center}
\renewcommand\arraystretch{1.2}
\caption{Physical parameters for Mandel’s problem.}
\label{PP}
\centering
{  \footnotesize
	\begin{tabular}{ccc}
	\hline
    Symbol & Description & Value   
    \\ \hline
    $a$       &  Dimension in $x$           &  $1.0$ $m$                             \\
    $b$       &  Dimension in $y$           &  $1.0$ $m$                             \\
    $h$       &  mesh size                  &  $0.1$ $m$                             \\
    $T$       &  Total simulation time      &  $1.0$ $s$                             \\
   $\Delta t$ &  Time step size             &  $0.001$ $s$                           \\
    $\lambda$ &  First Lam{\'e} parameter   &  $1.65 \times 10^9$ $Pa$               \\
    $\mu$     &  Second Lam{\'e} parameter  &  $2.475 \times 10^9$ $Pa$              \\
    $\alpha$  &  Biot-Willis constant       &  $1.0$                                 \\
    $c_0$     &  Specific storage           &  $6.061 \times 10^{-11}$ $Pa^{-1}$     \\
    $k_p$       &  Hydraulic conductivity     &  $9.869 \times 10^{-11}$ $m^2 Pa^{-1}$ \\
    $F$       &  Applied load               &  $6.0 \times 10^8$ $N m^{-1}$          \\
    $B$       &  Skempton's coefficient     &  $0.833$                               \\
    $\nu_u$   &  undrained Poisson’s ratio  &  $0.44$                                \\
    \hline
	\end{tabular}
}
\end{center}
\end{table}

\begin{figure}
\centering
\includegraphics[width=0.95\textwidth]{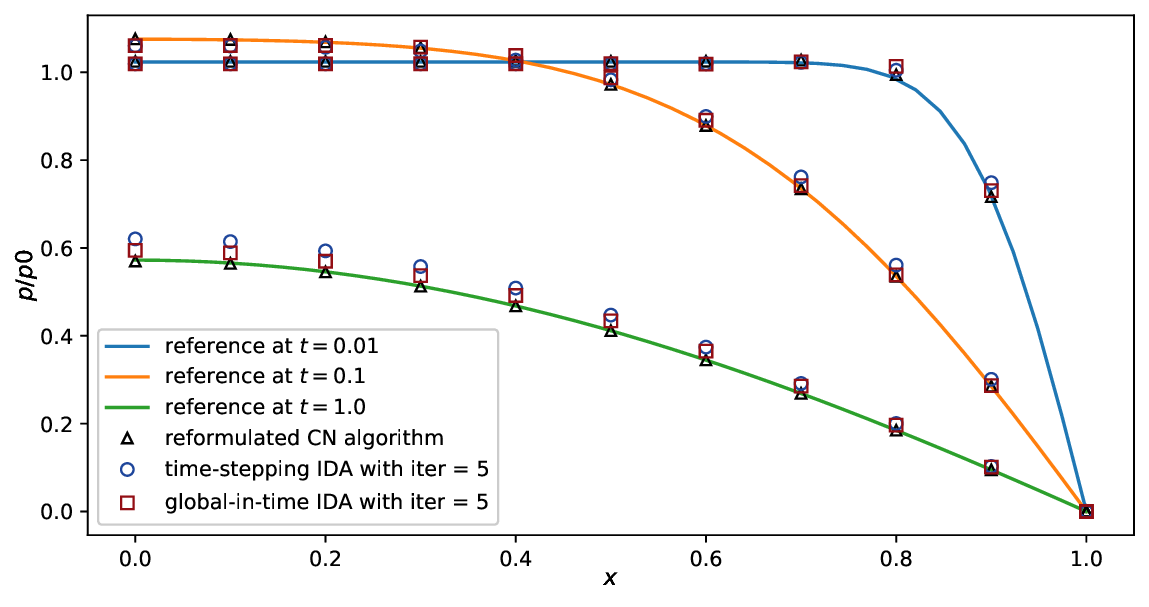}
\caption{Pressure solutions obtained by the coupled algorithm (\textbf{Algorithm \ref{algo:Coupled_CN}}), the proposed iterative decoupled algorithms (\textbf{Algorithm \ref{algo:TIDA_CN}} and \textbf{Algorithm \ref{algo:GIDA_CN}}) on Mandel’s problem.}
\label{figm}
\end{figure}

In Figure \ref{figm}, we compare the numerical results of the pressure obtained by Algorithms \ref{algo:Coupled_CN}, \ref{algo:TIDA_CN}, \ref{algo:GIDA_CN} with the analytical solution.
Here, we maintain the maximum iteration number $I = 5$ for the two iterative decoupled algorithms.
From Figure \ref{figm}, it is evident that \textbf{Algorithm \ref{algo:GIDA_CN}} performs slightly better than \textbf{Algorithm \ref{algo:TIDA_CN}}. 
This suggests that the global-in-time iterative decoupled algorithm may exhibit superior performance for solving long-time problems.
Furthermore, the results of \textbf{Algorithm \ref{algo:TIDA_CN}} and \textbf{Algorithm \ref{algo:GIDA_CN}} converge to the solution of the coupled algorithm as expected, which verifies our theoretical prediction.

\section{Conclusions}
\label{conclusion}

In this study, we integrated the traditional Crank-Nicolson method into the time discretization process of the three-field Biot's model, resulting in the development of the reformulated Crank-Nicolson algorithm. This reformulated algorithm inherits the advantageous features of unconditional stability and second-order time accuracy. Through the integration of the iterative decoupled algorithm with the reformulated Crank-Nicolson method, two convergent iterative decoupled algorithms emerged: the time-stepping Crank-Nicolson iterative decoupled algorithm and the global-in-time Crank-Nicolson iterative decoupled algorithm. It is shown that the solutions of time-stepping and global-in-time iterative decoupled algorithms converge to the solution of the coupled algorithm. Their effectiveness is verified through numerical experiments. We note that the global-in-time algorithm operates by solving solutions over the entire time interval within each iteration, which is a parallel-in-time computing strategy and is suitable for long-time problems. Notably, these methods, characterized by their simplicity and outstanding results, prove to be easily applicable in practical scenarios. 


\section{Acknowledgements}
The work of H. Gu is supported by the National NSF of China No. 123B2016. The work of J. Li is supported by the National NSF of China No. 11971221 and the Shenzhen Sci-Tech Fund No. RCJC20200714114556020, JCYJ20170818153840322 and JCYJ20190809150413261, Guangdong Provincial Key Laboratory of Computational Science and Material Design No. 2019B030301001, and GuangDong Basic and Applied Basic Research Foundation 2023B1515250005. The work of M. Cai is supported in part by an NIH-RCMI grant through U54MD013376,  the affiliated project award from the Center for Equitable Artificial Intelligence and Machine Learning Systems (CEAMLS) at Morgan State University (project ID 02232301). 


\appendix

\bibliographystyle{amsplain}
\bibliography{sample}

\end{document}